\documentclass[12pt]{amsart}
\usepackage{amssymb}
\ifx\mathrm\undefined\let\mathrm\rm\fi
\ifx\mathbf\undefined\let\mathbf\bf\fi
\ifx\mathfrak\undefined\let\mathfrak\frak\fi
\ifx\mathcal\undefined\let\mathcal\cal\fi
\ifx\mathbb\undefined\let\mathbb\Bbb\fi
\ifx\emph\undefined\let\emph\it\fi
\newcommand{\g}{{{\mathfrak g}\,}}

\newcommand{\n}{{{\mathfrak n}}}
\newcommand{\h}{{{\mathfrak h\,}}}
\newcommand{\OO}{{\mathcal O}}
\newcommand{\W}{{{\Bbb W\,}}}

\newcommand{\A}{{\mathcal A}}
\newcommand{\B}{{\mathcal B}}

\newcommand{\Z}{{\mathbb Z}}

\newcommand{\C}{{\mathbb C}}

\newcommand{\be}{\begin{displaymath}}
\newcommand{\ee}{\end{displaymath}}
\newcommand{\bea}{\begin{eqnarray*}}
\newcommand{\eea}{\end{eqnarray*}}

 \let\eps\varepsilon \let\epsilon\eps

\let\al\alpha

\let\ka\kappa
\let \la \lambda 
 
 \let\phi\varphi

\newcommand{\ovo}{{\overrightarrow{\otimes}}}

\newcommand{\bean}{\begin{eqnarray}}
\newcommand{\eean}{\end{eqnarray}}

\newcommand{\RR}{{\Bbb R}}

\newcommand{\T}{\!\otimes\!\,}

\newcommand{\vs}{\vspace{.5\baselineskip}}
\newenvironment{definition}
{\noindent{\bf Definition\/}:}{\par\vs}
\newtheorem%
{thm}{Theorem}%[section]
\newtheorem%
{proposition}[thm]{Proposition}
\newtheorem%
{lemma}[thm]{Lemma}
\newtheorem%
{lemmadef}[thm]{Lemma-Definition}
\newtheorem%
{corollary}[thm]{Corollary}
\newtheorem%
{conjecture}[thm]{Conjecture}
\title[ Dynamical Weyl groups and applications]
{Dynamical Weyl groups and applications}

\author[\smash{P.\;Etingof and A.\;Varchenko}]
{P.\;Etingof$^{\,\star,1}$ and
A.\;Varchenko$^{\,\diamond,2}$}
\thanks{$^1$ Supported in part by NSF grant DMS-9700477;
this research was partially conducted by the first
author for the Clay Mathematics Institute}
\thanks{$^2$ Supported in part by NSF grant DMS-9801582}

\begin{document}
\maketitle

\begin{center}
{\it
$^\star$   Department of Mathematics, MIT, Cambridge, MA
02139, USA,

and Columbia University,
Department of Mathematics,
2990 Broadway, New York, NY 10027, USA

etingof@math.mit.edu

\medskip

$^\diamond$Department of Mathematics, University of North Carolina
at Chapel Hill,

Chapel Hill, NC 27599-3250, USA,

av@math.unc.edu}
\end{center}
%\vsk1.8>

\begin{abstract}
Following a preceding paper of Tarasov
and the second author, we define and study
a new structure, which may be regarded as
the dynamical analogue of the Weyl group for Lie algebras
and of the quantum Weyl group for quantized enveloping algebras.
We give some applications of this new structure.
\end{abstract}

\thispagestyle{empty}

\section{Introduction}

In 1994, G.Felder, in his pioneering work \cite{F}, initiated
the development of a new area of the theory of quantum groups --
the theory of dynamical quantum groups. This theory assigns
dynamical analogs to various objects related to ordinary Lie
algebras and quantum groups (e.g. Hopf algebras, R-matrices, twists, etc.)
In particular, the main goal of
the present paper is to assign a dynamical analog to the Weyl
group of a Kac-Moody Lie algebra $\g$ and the quantum Weyl group
of the corresponding quantum group.
More specifically, we give a (rather straightforward)
generalization of the main construction of
the paper \cite{TV}, which, in effect, introduces
dynamical Weyl groups in the case of finite dimensional
simple Lie algebras.

The analog of the Weyl group we introduce
is a collection of operators
that give rise to a braid group representation on the space of functions
from the dualized Cartan subalgebra $\h^*$ of $\g$ to a representation
$V$ of $\g$ or $U_q(\g)$.
We call this analog
{\bf the dynamical Weyl group} of $V$.

We note that dynamical Weyl groups may be regarded as generalizations
of the classical ``extremal projectors'' introduced in \cite{AST}.
In particular, the dynamical Weyl group operators for simple Lie algebras
were introduced in \cite{Zh1,Zh2}, by analogy with \cite{AST}
(see formula (3.5) and
Theorem 2 in \cite{Zh2}). This construction, however, is
different from
that of \cite{TV}.

Dynamical Weyl groups are not only beautiful
objects by themselves, but also have a number of useful
applications. To describe one of these applications, recall that in
\cite{EV2}, we developed the theory of
trace functions (matrix analogs of Macdonald functions),
using the basic dynamical
objects introduced in \cite{EV1} (the fusion and exchange
matrices). In particular, we derived
four systems of difference equations for these functions:
qKZB, dual qKZB, Macdonald-Ruijsenaars, and dual
Macdonald-Ruijsenaars equations,
and proved the symmetry of the trace functions
under permutation of components and arguments simultaneously.
In this paper, we use the dynamical Weyl group
to develop this theory further: namely, we show that
the trace functions and all four systems of equations for them
are symmetric with respect to the dynamical Weyl group
(while they are not, in general, symmetric under the usual
classical or quantum Weyl group). This property is a
generalization to the matrix case of the Weyl symmetry
property of Macdonald functions; being important by itself, it
also allows one to prove other properties of trace functions
(orthogonality, the Cherednik-Macondald-Mehta identities),
which we plan to do in a separate paper.

As a second application, we interpret the
important operator $Q(\la)$ from \cite{EV2} in terms of
the dynamical Weyl group operator corresponding to the
maximal element of the Weyl group. This allows us to calculate
$Q(\la)$ explicitly, and in particular get an explicit
product formula for its determinant. In the next paper, we
will show that $Q(\la)$ is the matrix analog of the
squared norm of the Macdonald polynomial $P_\la$,
and in particular the product formula for the determinant of
$Q(\la)$ is the matrix analog of the well known
Macdonald inner product
identities.

Finally, a third application of dynamical Weyl groups is to the theory of
KZ and qKZ equations, and is along the lines of \cite{TV}.
Recall that the main goal of \cite{TV} was not to define
dynamical Weyl groups for simple Lie algebras, but rather
to construct commuting difference operators which commute with
the trigonometric KZ operators and are a deformation of the
differential operators from \cite{FMTV} commuting with the KZ
operators. Such difference operators were constructed in
\cite{TV} using the dynamical Weyl group of the
corresponding simple Lie algebra, combined with the method of
Cherednik \cite{Ch1} of lifting R-matrices to affine R-matrices.
However, the method of \cite{TV} allowed the authors of \cite{TV} to prove
that their operators actually commute with the trigonometric KZ operators
only for Lie algebras other than $E_8,F_4,G_2$
(i.e. for Lie algebras having a minuscule fundamental coweight).
In this paper, we attack the same problem using a somewhat
different method: we use the dynamical Weyl group of the {\it affine}
Lie algebra, rather than the finite dimensional one, which allows us
to avoid using the procedure from \cite{Ch1}.
As a result, we obtain the same difference operators as in \cite{TV},
and prove that they commute with the KZ operators 
for any simple Lie algebra $\g$. 
We also generalize the constructions and results
of \cite{TV} to the quantum group case.

We note that another unexpected application of dynamcial Weyl
groups appears in the recent interesting paper \cite{RS}.

The organization of the paper is as follows.

In Section 2, we recall the basic notions used in this
paper: Kac-Moody algebras and their quantizations,
(ordinary) Weyl groups, intertwining operators, Verma modules,
singular vectors, fusion and exchange matrices.

In Section 3, using the operation of restriction of intertwining operators
for Verma modules to their submodules, we define the dynamical Weyl group
operator for $U_q({{\frak sl}_2})$. We calculate this operator
explicitly.

In Section 4, we use the operator from Section 3 to define
the dynamicl Weyl group operators for any quantized
Kac-Moody algebra. We show that similarly to the $U_q({{\frak sl}_2})$ case,
the dynamical Weyl group arises from the restriction procedure
 for intertwiners.

In Section 5, we continue to study the properties of
the dynamical Weyl group, and, in particular, show that its limit
at infinity gives the usual quantum Weyl group of Soibelman and
Lusztig.

In Section 6 we give the first application of the dynamical Weyl
group: we link the operator $Q(\lambda)$ from \cite{EV2}
with the dynamical Weyl group operator corresponding to the
maximal element of the Weyl group.

In Section 7, we describe the applications of the
dynamical Weyl group to trace functions. We establish the
dynamical Weyl group symmetry of these functions and of the
equations for them introduced in \cite{EV2}.

In Section 8, we discuss the dynamical Weyl groups
of loop representations
for affine Lie algebras and quantum affine algebras.

In Section 9, using the material of Section 8,
we give a more conceptual derivation of
the difference equations from \cite{TV}, which commute
with the trigonometric KZ equations. We also derive the q-analogs
of the equations from \cite{TV}, which commute with the
trigonometric quantum KZ equations.

{\bf Acknowledgements.}
We are grateful to S. Khoroshkin, A. Kirillov Jr., G. Lusztig, Ph. Roche,
and V.Tarasov for useful discussions and references. The first author thanks
IHES for hospitality.

\section{Preliminaries}

\subsection{Kac-Moody algebras}

We recall definitions from \cite{K}.
Let $A=(a_{ij})$ be a
symmetrizable generalized Cartan matrix of size $r$,
and $(\h,\Pi,\Pi^\vee)$ be a realization of $A$.
This means that
$\h$ is a vector space of dimension $2r-\text{rank}(A)$,
$\Pi=\{\alpha_1,...,\alpha_r\}\subset \h^*$,
$\Pi^\vee=\{h_1,...,h_r\}\subset \h$ are
linearly independent, and $\alpha_i(h_j)=a_{ji}$.
The elements $\alpha_i$ are called simple positive roots.

\begin{definition} The Kac-Moody Lie algebra $\g(A)$ is generated
by $\h,e_1,...,e_r,f_1,...,f_r$ with defining relations
$$
[h,h']=0,\ h,h'\in\h;\ [h,e_i]=\alpha_i(h)e_i; \
[h,f_i]=-\alpha_i(h)f_i;\ [e_i,f_j]=\delta_{ij}h_i,
$$
and the Serre relations
$$
\sum_{m=0}^{1-a_{ij}}\frac{(-1)^m}
{m!(1-a_{ij}-m)!}e_i^{1-a_{ij}-m}e_je_i^m=0,
$$
$$
\sum_{m=0}^{1-a_{ij}}\frac{(-1)^m}
{m!(1-a_{ij}-m)!}f_i^{1-a_{ij}-m}f_jf_i^m=0.
$$
\end{definition}

For brevity we will assume that $A$ is fixed and denote $\g(A)$
simply by $\g$. The positive and negative nilpotent subalgebras
of $\g$ will be denoted by $\n_\pm$.

By the definition of a
generalized Cartan matrix, there exists a collection of positive integers
$d_i$, $i=1,...,r$,
such that $d_ia_{ij}=d_ja_{ji}$. We will choose the minimal collection of
such numbers, i.e. the collection for which the numbers are
the smallest possible (such a choice is unique).
Let us choose a nondegenerate
bilinear symmetric form on $\h$ such that
$(h,h_i)=d_i^{-1}\alpha_i(h)$.
It is easy to see that such a form always exists.
It is known \cite{K} that there exists a unique extension of the
form $(,)$ to an invariant nondegenerate symmetric bilinear form
$(,)$ on
$ \g$. For this extension, one has
$(e_i,f_j)=\delta_{ij}d_i^{-1}$.

{\bf Remark.}
One can show that forms on $\g$ coming from different
forms on $\h$ are equivalent under automorphisms of $\g$.

A root of $\g$ is a nonzero element of $\h^*$ which occurs in the
decomposition of $\g$ as an $\h$-module. A root is positive if it
is a positive linear combination of simple positive roots, and
negative otherwise. A root $\alpha$ is real if
$(\alpha,\alpha)>0$,
otherwise it is imaginary.
For a real root $\alpha$ of $\g$, let $\alpha^\vee=
2\alpha/(\alpha,\alpha)$ be the corresponding coroot.

\subsection{Quantized Kac-Moody algebras}

Let
$\hbar$ be a complex number, which is not
a rational multiple of $\pi i$,
and $q=e^{\hbar/2}$. For a number or operator $B$, by
$q^B$ we mean $e^{\hbar B/2}$.

\begin{definition} The quantized Kac-Moody algebra
$U_q(\g(A))$ is
the associative algebra generated
by $e_1,...,e_r,f_1,...,f_r$, and $q^h,h\in\h$
(where $q^0=1$), with defining relations
$$
q^hq^{h'}=q^{h+h'},\ h,h'\in\h;\ q^he_i=q^{\alpha_i(h)}e_iq^h; \
q^hf_i=q^{-\alpha_i(h)}f_iq^h;\ [e_i,f_j]=\delta_{ij}
\frac{q_i^{h_i}-q_i^{-h_i}}{q_i-q_i^{-1}},
$$
and the Serre relations
$$
\sum_{m=0}^{1-a_{ij}}\frac{(-1)^m}
{[m]_{q_i}![1-a_{ij}-m]_{q_i}!}e_i^{1-a_{ij}-m}e_je_i^m=0.
$$
$$
\sum_{m=0}^{1-a_{ij}}\frac{(-1)^m}
{[m]_{q_i}![1-a_{ij}-m]_{q_i}!}f_i^{1-a_{ij}-m}f_jf_i^m=0,
$$
where $[m]_q=\frac{q^m-q^{-m}}{q-q^{-1}}$,
and $q_i:=q^{d_i}$.
\end{definition}

For brevity we will denote $U_q(\g(A))$
by $U_q(\g)$. Also, to give a uniform treatment
of the classical and quantum case, we will
often allow $\hbar$ to be $0$ (i.e. $q=1$),
in
which case $U_q(\g)$ is defined to be $U(\g)$.

The positive and negative nilpotent subalgebras in $U_q(\g)$ will
be denoted by $U_q(\n_\pm)$.

The algebra $U_q(\g)$ is a Hopf algebra, with coproduct defined
by
$$
\Delta(q^h)=q^h\otimes q^h, \Delta(e_i)=e_i\otimes
q_i^{h_i}+1\otimes e_i, \Delta(f_i)=f_i\otimes
1+q_i^{-h_i}\otimes f_i,
$$
and the antipode defined by
$$
S(e_i)=-e_iq_i^{-h_i},\ S(f_i)=-q_i^{h_i}f_i,\ S(q^h)=q^{-h}.
$$

\subsection{Verma modules and integrable modules}

Let $\lambda\in \h^*$ be a weight.
We say that a vector $v$ in a module $V$ over $\g$ or $U_q(\g)$
has weight $\lambda$ if $hv=\lambda(h)v$ for all $h\in \h$
(respectively $q^hv=q^{\lambda(h)}v$).
The space of vectors of weight $\lambda$ is denoted by
$V[\lambda]$. Modules in which any
vector is a sum of vectors of some weight are said to be
$\h$-diagonalizable. Category $\OO$ consists of $\h$-diagonalizable
modules with finite dimensional weight subspaces, whose weights
belong to a union of finitely many ``conical'' sets of the form
$\lambda-\sum_i\Z_+\alpha_i$.

An example of an object in $\OO$ is a Verma module.
The Verma module
$M_\lambda$ over $\g$
with highest weight $\lambda$ is generated by one generator
$v_\la$ with
defining relations
$e_iv_\la=0$, $hv_\la=\lambda(h)v_\la$, $h\in \h$.
The Verma module $M_\lambda$ over $U_q(\g)$ is generated by $v_\la$
with
defining relations $e_iv_\la=0$, $q^hv_\la=q^{\lambda(h)}v_\la$.

{\bf Remark.} All the Verma modules
in this paper are equipped with a distinguished generator
(i.e. the normalization of the generator is fixed).

We say that an object $V$ in $\OO$ is integrable if for all $i$, it is a sum
of finite dimensional submodules with respect to the subalgebra
generated by $e_i,f_i,q^{bh_i}$, $b\in \C$.

We say that $\lambda$ is a dominant integral weight if
$\lambda(h_i)$ is a nonnegative integer for all $i$.
The set of dominant integral weights is denoted by $P_+$.

The irreducible module $L_\lambda$ with highest
weight $\lambda$ over $U_q(\g)$ is an integrable module if and only if
$\lambda$ is a dominant integral weight.
The category $\OO_{int}$ of integrable modules is semisimple, with
irreducible objects being $L_\lambda$ for dominant integral
$\lambda$. This category is closed under tensor product.

\subsection{The Weyl group}

Recall that the Weyl group $\W$ of $\g$ is the group of
transformations of $\h$ generated by
the reflections $s_i(\lambda)=\lambda-\lambda(h_i)\alpha_i$.
It is known that the defining relations for $\W$ are: $s_i^2=1$,
$(s_is_j)^{m_{ij}}=1$, $i\ne j$, where $m_{ij}=2,3,4,6,\infty$ if
$a_{ij}a_{ji}=0,1,2,3,\ge 4$ (if $m_{ij}=\infty$ then we agree
that there is no relation).
Any element of $\W$ is a product of $s_i$.
The smallest number of factors in such a product is called the
length of $w$ and denoted by $l(w)$. A representation of
$w\in \W$ by a product of length $l(w)$ is called a reduced
decomposition.

The group $\W$ is the quotient of the group $\tilde\W$ generated
by $s_i$ with the braid relations
$$
s_is_js_i...=s_js_is_j...
$$
($m_{ij}$ terms on both sides), by the additional relations
$s_i^2=1$. The group $\tilde \W$ is called the braid group of
$\g$.

Any two reduced decompositions of an element
of $\W$ coincide not only in $\W$ but also in $\tilde \W$
(see e.g. \cite{Lu}, 2.1.2). Therefore,
the projection map $\tilde \W\to \W$ admits a
splitting $\gamma: \W\to \tilde \W$, assigning to any element
of $\W$ given by some reduced decomposition, the element
of $\tilde \W$ defined by the same decomposition
(of course, this is only a map of sets, not a group homomorphism).
Using the map
$\gamma$, we will regard $\W$ as a subset of $\tilde \W$.

Let us fix a weight $\rho$ such that $\rho(h_i)=1$ for all $i$.
Define the shifted action of the Weyl group on $\h$ by $w\cdot
\lambda=w(\lambda+\rho)-\rho$. It is obvious that the shifted
action is independent on the choice of $\rho$.

\subsection{Intertwining operators and expectation values}

Let $V$ be an $\h$-diagonalizable module over $U_q(\g)$,
and $\Phi: M_\lambda\to M_\mu\otimes V$ an intertwining
operator. We have $\Phi v_\lambda=v_\mu\otimes v+...$, where $...$
denote terms of lower weight in the first component, and
$v\in V[\lambda-\mu]$. We will call $v$ the expectation value of
$\Phi$ and write $<\Phi>=v$.

Let $V$ be a module over
$U_q(\g)$ which belongs to $\OO$. Let $\nu$ be a weight of $V$.

\begin{lemma}\label{inter}
For generic $\la$ the map $\Phi\to <\Phi>$ is an isomorphism of
  vector spaces
$\text{Hom}_{U_q(\g)}(M_\lambda,M_{\lambda-\nu}\otimes V)\to
V[\nu]$. In particular, this
map is an isomorphism for dominant weights
with sufficiently large coordinates $\la(h_i)$ for all
$i=1,...,r$.
\end{lemma}

\begin{proof} The proof is straightforward; see \cite{ES} and
\cite{ESt}.
\end{proof}

This lemma allows one to define, for $v\in V[\mu]$ and
generic $\lambda\in \h^*$,
 the intertwining operator
$\Phi_\lambda^v$ such that $<\Phi_\lambda^v>=v$.
It is easy to see that the matrix elements of this operator
with respect to the bases in $M_\lambda$,
$M_{\lambda-\mu}$ induced by any bases in $U_q(\n_-)$ and $V$,
are rational functions of $(\lambda,\alpha_i)$ for $q=1$ and of
$q^{(\lambda,\alpha_i)}$ if $q\ne 1$.

\subsection{Fusion and exchange matrices}

Recall the definition of the fusion and exchange matrices
\cite{ES,EV1}.

Let $\la\in \h^*$ be a generic weight.
Let $V,U$ be integrable $U_q(\g)$-modules, and
$v \in V[\mu],\;u\in U[\nu]$.
Consider the composition
\bean
\Phi^{u,v}_{\lambda}:\;M_\lambda
\stackrel{\Phi^v_\lambda\otimes 1}{\longrightarrow} M_{\lambda-\mu}
\otimes V \stackrel{\Phi^u_{\lambda-\mu}\otimes 1}{\longrightarrow}
M_{\lambda-\mu-\nu} \otimes U \otimes V.
\notag
\eean
Then
$\Phi^{u,v}_\lambda \in
\mathrm{Hom}_{U_q(\g)}(M_\lambda,M_{\lambda-\mu-\nu}
\otimes U \otimes V)$.  We will call $\Phi_\la^{u,v}$ the {\it fusion} of
$\Phi_\la^v$ and $\Phi_{\la-\mu}^u$.

For a generic $\lambda$
there exists a unique element $x \in \,(U \otimes
V) [\mu+\nu]$ such that
$\Phi^x_\lambda=\Phi^{u,v}_\lambda$. The assignment
$(u, v) \mapsto x$ is bilinear, and defines a zero weight map
$$
J_{UV}(\lambda):\; U \otimes V \to U \otimes V.
$$
The operator $J_{UV}(\lambda)$ is called the fusion matrix of $U$ and $V$.
The fusion matrix $J_{UV}(\lambda)$ is a
rational function of $\lambda$
for $q=1$ (respectively of $q^\lambda$ if $q\ne 1$).

Also, $J_{UV}(\lambda)$ is strictly
lower triangular, i.e. $J=1+L$ where $ L(U[\nu] \otimes V[\mu])
\subset \oplus_{\tau<\nu, \,\mu<\sigma} U[\tau]\otimes V[\sigma]$.
In particular, $J_{UV}(\lambda)$ is invertible.

The exchange matrix $R_{VU}(\lambda)$ is defined by the formula
$$
R_{VU}(\lambda)=J_{VU}(\lambda)^{-1}{\mathcal
  R}^{21}J^{21}_{UV}(\lambda),
$$
where ${\mathcal R}$ is the universal R-matrix of $U_q(\g)$.
The exchange matrix has zero weight. It is also called a
dynamical R-matrix, since it
satisfies the quantum dynamical Yang-Baxter
equation (see \cite{EV1}).

In the theory of fusion and exchange matrices, one often
uses the so called ``dynamical notation'', which will also
be useful for us here. This notation is defined as follows.

Let $V_1,\ldots V_n$ be $\h^*$-graded vector spaces, and let $F(\lambda):
V_1 \otimes \ldots \otimes V_n \to V_1 \otimes \ldots \otimes V_n$ be a
linear operator depending on $\lambda \in \h^*$. Then for any
homogeneous $u_1,\ldots , u_n$, $u_i\in V_i[\nu_i]$,
 we  define $ F(\lambda-h^{(i)})(u_1 \otimes
\ldots \otimes u_n)$ to be $F(\lambda-\nu_i) (u_1 \otimes \ldots
\otimes u_n)$. In particular, when $n=1$, we will denote the term
$h^{(1)}$ simply by $h$: that is, $F(\la-h)v=F(\la-\nu)v$ if
$v$ has weight $\nu$.

\subsection{Singular vectors in Verma modules.}

Recall that a nonzero vector in a $U_q(\g)$-module
is said to be singular if it is annihilated by $e_i$ for all $i$.

Let $w\in \W$ and $w=s_{i_1}\ldots s_{i_l}$ be a reduced
decomposition. Set $\al^{l}=\al_{i_l}$ and $\al^{j}=(s_{i_l}\ldots
s_{i_{j+1}}) (\al_{i_j})$ for $j=1,\ldots,l-1$. Let
$n_j=2\frac{(\la+\rho,\alpha^j)} {(\alpha^j,\alpha^j)}$. For a
dominant $\la\in P_+$, $n_j$ are positive integers. Let
$d^j=d_{i_j}$ (where $d_i$ are the symmetrizing numbers).

We will need the following lemma, which is similar to
Lemma 4 in \cite{TV}.

\begin{lemma}\label{sing-v} Let $\la$ be a
dominant integral weight. Then
the collection of pairs of integers $(n_1,d^1),..., (n_k,d^k)$ and
the product $f_{\al_{i_1}}^{n_1}\cdots f_{\al_{i_l}}^{n_l}$ does
 not depend on
the reduced decomposition.
\end{lemma}

{\bf Remark.} The second statement of the lemma
is known as the ``quantum Verma identities''
and is proved in \cite{Lu}, Section 39.3.
In the case $q=1$, it goes back to Verma.
However, we will give the
argument (which is somewhat different from the one in \cite{Lu})
for the reader's convenience.

\begin{proof} It is sufficient to prove the statement
for two reduced decompositions that can be identified by applying
a braid relation once. Therefore, it is sufficient to check
the statement for rank 2 Lie algebras $A_2,B_2,G_2$. In this case
the only element that has two different reduced decompositions is
the maximal element $w_0$. So we can assume that we are dealing
with the two different reduced decompositions of $w=w_0$, namely,
$w=s_{i_1}...s_{i_l}=s_{i_1'}...s_{i_l'}$.

Since $w=w_0$, it is easy to see that for either reduced
decomposition, $\al^{1},\ldots,\al^{l}$ are all the positive roots
(each occuring exactly once).

Hence, the collection $(n_1,d^1)\ldots (n_l,d^l)$ does not depend
on the reduced decomposition.

Let $n_i,n_i'$ be the numbers defined above for the two
decompositions. The vectors $u=f_{\al_{i_1}}^{n_1}\cdots
f_{\al_{i_l}}^{n_l}v_\la$, $u'=f_{\al_{i_1'}}^{n_1'}\cdots
f_{\al_{i_l'}}^{n_l'}v_\la$ are singular vectors in $M_\la$ of
weight $w_0\cdot \la$ (they are nonzero, since the algebra
$U_q(\n_-)$ has no zero divisors).

For Lie algebras $A_2,B_2,G_2$, it is easy to see
that the space of singular vectors in $M_\la$ in the weight $w_0\cdot \la$
is 1-dimensional.
Indeed, the module
$M_{w_0\cdot \la}$ is irreducible, so if there
were two independent singular vectors of weight $w_0\cdot \la$,
then the direct sum of two copies of $M_{w_0\cdot \la}$
would be contained in $M_\la$. But this is impossible,
since some weight multiplicities of this direct sum are bigger
than the corresponding weight multiplicities of $M_\la$.

Therefore, the vectors $u,u'$ are proportional. Since $M_\la$ is
a free module over the subalgebra $U_q(\n_-)$, we have
$f_{\al_{i'_1}}^{n'_1}\ldots f_{\al_{i'_l}}^{n'_l}\,=\,
c\,f_{\al_{i_1}}^{n_1}\ldots f_{\al_{i_l}}^{n_l}$ in $U_q(\n_-)$
for a suitable $c\in\C^*$.

We claim that $c=1$. Indeed, consider the
 natural homomorphism from
$U_q(\n_-)$ to the algebra generated by $f_i$ with the relations
$f_if_j=q_i^{a_{ij}}f_jf_i$, $i<j$, sending $f_i$ to $f_i$ (it
is easy to check that such a homomorphism exists).
The images of the two monomials
under this homomorphism differ by a power of $q$,
which implies that $c=q^m$. On the other hand, since the Serre relations
are symmetric under $q\to q^{-1}$, a similar
 homomorphism exists if $q$ is replaced with $q^{-1}$, which
 yields $c=q^{-m}$. Thus, $c=1$, as desired.
\end{proof}

Let $\delta$ be a reduced decomposition of $w\in \W$ given by the
formula $w=s_{i_1}...s_{i_l}$. Define a vector
$v_{w\cdot\la,\delta}^\la\in M_\la$ by \bean
v_{w\cdot\la,\delta}^\la\,=\, \frac{f_{\al_{i_1}}^{n_1}}{
[n_1]_{q^{d^1}}!} \ldots \frac{f_{\al_{i_l}}^{n_l}}{
[n_l]_{q^{d^l}}!}\,v_\la\,, \eean This vector is singular. It does
not depend on the reduced decomposition $\delta$ by Lemma
\ref{sing-v}, so we will often denote it by $v_{w\cdot \la}^\la$.

\section{The main construction for $U_q({{\frak sl}_2}$)}

\subsection{The operators $A_{s,V}(\la)$}

Let $\g={{\frak sl}_2}$.
Identify the space of weights for ${{\frak sl}_2}$ with $\C$ by $z\in
\C\to z\alpha/2$, where $\alpha$ is the positive root;
then dominant integral weights are identified with positive integers.

Let $s$ be the nontrivial element of the Weyl group of ${{\frak sl}_2}$.
Let $V$ be a finite dimensional $U_q({{\frak sl}_2})$-module,
and let $\la$ be a sufficiently large positive integer
(compared to $V$).
 Define
a linear operator $A_{s,V}(\lambda): V\to V$ as
follows.

Fix a weight $\nu$ of $V$. Let $v\in V[\nu]$.
Consider the intertwining operator
$\Phi^v_\la: M_{\la}\to M_{\la-\nu}\T V$.
By Lemma \ref{inter}, such an operator is well defined.

\begin{lemma}\label{a} For a sufficiently large positive integer
$\la$,
there exists a unique linear operator $A_{s,V}(\la):V\to V$ such that
\bean\label{A}
\Phi_\la^v v_{s\cdot\la}^\la= v_{s\cdot(\la-\nu)}^{ \la-\nu}\T
 A_{s,V}(\la)v\,
+\,\text{lower weight terms}\,.
\notag
\eean
(where the weight is taken in the first component).
This operator is invertible.
\end{lemma}

\begin{proof} The proof of existence and uniqueness
of $A_{s,V}$ is
straightforward (see e.g. \cite{TV}).
To prove the invertibility,
it is sufficient to observe
that for an irreducible module $V$ and large
$\la$, the map $A_{s,V}(\la):V[\nu]\to V[-\nu]$
is nonzero. Indeed, the tensor product of a Verma module with a finite
dimensional module does not contain finite dimensional
submodules. Thus, the operator $\Phi_\la^v$ cannot have finite rank,
and hence has to be nonzero on $M_{-\la-2}$.
\end{proof}

Thus, the operator $A_{w,V}(\la)$ is
 the effect, at the level of expectation values,
of the operation of {\bf restriction} of an intertwiner from
$M_\la$ to $M_{s\cdot \la}$.

The goal of the next few sections is to compute the operator
$A_{s,V}(\lambda)$. This can be done
analogously to \cite{TV}, using a direct
calculation and identities with hypergeometric functions.
However, we would like
to give a different derivation, which seems to be a bit
simpler (in the spirit of \cite{EV2}, subsection 7.2).

The main tool of the calculation is
the following important property of $A_{s,V}(\la)$.

\begin{lemma}\label{com-Asl2}
Let $U, V$ be finite dimensional
$U_q({{\frak sl}_2})$-modules. Then
\bean\label{Com-Asl2}
A_{s,U\T V}(\la)J_{UV}(\la)\,=\,J_{UV}(s\cdot\la)
A_{s,V}^{(2)}(\la)A_{s,U}^{(1)}(\la-h^{(2)})
\,
\eean
where $A^{(1)}$ denotes $A\otimes 1$,
$A^{(2)}$ denotes $1\otimes A$.
\end{lemma}

\begin{proof} The lemma is an easy consequence of the
  definitions: it expresses the fact that
the operation of fusion of intertwiners
  commutes with the operation of restriction
of intertwiners to submodules.
\end{proof}

\subsection{Calculation of $A_{s,V}(\la)$ in the 2-dimensional
  representation}

For brevity we denote $A_{s,V}$ by $A_V$. Consider the case when
$V$ is the 2-dimensional irreducible representation
with the standard basis $v_+$ and $v_-$,
such that $ev_+=fv_-=0$, $ev_-=v_+$, $fv_+=v_-$,
$q^{bh}v_\pm=q^{\pm b}v_\pm$.

\begin{lemma} One has
$$
A_V(\la)v_+=qv_-,\quad A_V(\la)v_-=-q^{-1}\frac{[\la+2]_q}{[\la+1]_q}v_+.
$$
\end{lemma}

\begin{proof}
Consider the intertwiner $\Phi_\la^{v_+}$. It satisfies
the relation
$$
\Phi_\la^{v_+}v_\la=v_{\la-1}\otimes v_+.
$$
Therefore,
$$
\Phi_\la^{v_+}\frac{f^{\la+1}}{[\la+1]_q!}v_\la=
\frac{1}{[\la+1]_q!}(f\otimes 1+q^{-h}\otimes
f)^{\la+1}(v_{\la-1}
\otimes v_+)
$$
This implies after a straightforward calculation that
$$
A_V(\la)v_+=qv_-.
$$

Now let us consider
the intertwiner $\Phi_{\la}^{v_-}$.
It satisfies
the relation
$$
\Phi_\la^{v_-}v_\la=v_{\la+1}\otimes v_--q^{-1}[\la+1]_q^{-1}
fv_{\la+1}\otimes v_+.
$$
Therefore,
$$
\Phi_\la^{v_-}\frac{f^{\la+1}}{[\la+1]_q!}v_\la=
\frac{1}{[\la+1]_q!}(f\otimes 1+q^{-h}\otimes
f)^{\la+1}(v_{\la+1}
\otimes v_--q^{-1}[\la+1]_q^{-1}
fv_{\la+1}\otimes v_+).
$$
This implies that
$$
A_V(\la)v_-=-q^{-1}\frac{[\la+2]_q}{[\la+1]_q}v_+,
$$
as desired.
\end{proof}

\subsection{The calculation in any finite dimensional
  representation (up to a constant)}

Now we let $V=V_m$ be the irreducible representation
with highest weight $m$.
For any $k=0,...,m$, define a linear map
$A_m^k(\lambda):V[m-2k]\to V[2k-m]$
to be the restriction of $A_V(\lambda)$ to $V[m-2k]$.
If we choose generators of the 1-dimensional
spaces $V[m-2k]$, this map
will be expressed by a scalar complex valued function of $\lambda$.

Up to a $\lambda$-independent factor, this function is independent
on the choice of the generators.
Thus, we can naturally understand
$A_m^k(\lambda)$ as an element of the group
$$
(\text{nonvanishing complex valued functions on } \Bbb Z_++N)/\C^*.
$$
(where $N$ is a large enough number).
This will be our point of view in this subsection.
The equality of two elements in this group will be denoted by
the sign $\equiv$.

\begin{proposition}\label{uptoconst}
One has
$$
A_m^k(\la)\equiv\prod_{j=1}^k\frac{[\la+1+j]_q}{[\la-m+k+j]_q}.
$$
\end{proposition}

{\bf Remark.}
In the case $q=1$, this proposition
appears in \cite{TV}.

\begin{proof}
Let $m\ge 1$.
Consider equation (\ref{Com-Asl2}) in the weight subspace of
weight $m-2k+1$ in the tensor product $V_1\otimes V_m$. Let
us identify this weight subspace with the opposite one in any
way, and take the determinant of both sides of (\ref{Com-Asl2}).

Since the fusion matrix is triangular
with the diagonal elements equal to $1$, its
determinant is $1$. Therefore, using the decomposition
$V_1\otimes V_m=V_{m-1}\oplus V_{m+1}$ we obtain for $k=0$:
$$
A_{m+1}^0(\lambda)\equiv A_m^0(\lambda)A_1^0(\lambda-m),
$$
and for $k\ne 0$
$$
A_{m+1}^k(\lambda)A_{m-1}^{k-1}(\lambda)\equiv
A_m^k(\lambda)A_m^{k-1}(\lambda)A_1^0(\lambda-m+2k)A_1^1(\lambda-m+2k-2).
$$
Now let us substitute the values of $A_1^0$ and $A_1^1$ computed
in the previous section. Then we get
$$
A_m^0(\lambda)=1,
$$
$$
A_{m+1}^k(\lambda)\equiv \frac{[\la-m+2k]_q}
{[\la-m+2k-1]}_q A_m^k(\la)A_m^{k-1}(\la)A_{m-1}^{k-1}(\la)^{-1}
$$
It is clear that $A_m^k$ is completely determined from this
equation. It remains to check that the expression given in the
proposition satisfies the equations, which is straightforward.
\end{proof}

\begin{corollary}
The operator-valued function $A_{V}(\la)$,
defined for large positive integers, uniquely extends
to a rational function of $\lambda$
(for $q=1$) and of $q^{\lambda}$ (for $q\ne 1$).
For generic $\lambda$, the operator $A_{s,V}(\lambda)$ is
invertible.
\end{corollary}

\begin{proof} The existence follows
from the above explicit computation of $A_{V}$.
The uniqueness is obvious, since the function is
  defined at infinitely many points.
The invertibility follows from Lemma \ref{a}.
\end{proof}

\subsection{Limits of $A_V(\la)$ at infinity}

In the previous subsection, we have calculated $A_V$ up to a constant.
In this subsection, we will explicitly calculate this
constant.

Let $V$ be a finite dimensional representation of $U_q({{\frak sl}_2})$.
Proposition \ref{uptoconst}
implies the following result.

\begin{corollary} \label{limits}
(i) If $q=1$ then the map $A_V(\lambda)$
has a limit $A_V^\pm=A_V^\infty$ at $\lambda=\pm\infty$.
If $q\ne 1$, the map $A_V(\lambda)$ has a limit $A_V^+$ as
$q^{\la}\to \infty$ and $A_V^-$ as $q^\la\to 0$, respectively.

(ii) One has
$$
A_V^-|_{V_m[m-2k]}=q^{-2k(m-k+1)}
A_V^+|_{V_m[m-2k]}.
$$
In other words, one has
$$
A_V^-=A_V^+\bold uq^{-h^2/2},
$$
where
$\bold u$ is the Drinfeld element $m_{21}(S\otimes 1){\mathcal
  R}$,
\cite{Dr}
(here $m_{21}(a\otimes b)=ba$).

(iii) Define
$$
A_V^\infty|_{V_m[m-2k]}=q^{-k(m-k+1)}A_V^+|_{V_m[m-2k]}=
q^{k(m-k+1)}A_V^-|_{V_m[m-2k]}.
$$
to be the geometric mean of the two limits.
Then
one has
$$
A_V(\la)=A_V^\infty B_V(\la),
$$
where $B_V(\la)$ is a weight zero operator, defined by
$$
B_{V_m}(\la)|_{V_m[m-2k]}=
\prod_{j=1}^k\frac{[\la+1+j]_q}{[\la+1-m+2k-j]_q}.
$$
\end{corollary}

{\bf Remark.} Recall that the element $\bold u$ acts
on $V_m$ as $q^{-m(m+2)/2}q^h$.

\subsection{Computation of $A_V^\pm,A_V^\infty$}

Let
 ${\mathcal R}_0={\mathcal R}q^{-h\otimes h/2}$.

\begin{proposition}\label{tensprod}
One has
$$
A_{U\otimes V}^+={\mathcal R}_0^{21}(A_U^+\otimes
A_V^+).
$$
\end{proposition}

\begin{proof}
Theorem 50 from \cite{EV1} implies that
as $q^{\la}\to \infty$, one has
$J_{UV}(\la)\to 1$ and
$J_{UV}(-\la)\to {\mathcal R}_0^{21}$. Thus,
going to the limit $q^{\la}\to \infty$ in (\ref{Com-Asl2}),
and using Theorem 50 from \cite{EV1}, we obtain
the proposition.
\end{proof}

{\bf Remark 1.} We take this opportunity
to correct the statement of Theorem 50 of \cite{EV1}.
First of all, the condition $|q|<1$ should be replaced
with $|q|>1$. (The paper \cite{ESt},
whose results are used to prove the theorem, refers to
the comultiplication opposite to that of \cite{EV1},
and the two comultiplications are related by the transformation
$q\to q^{-1}$; cf. formula (45) in \cite{EV1}). Second,
$\n_\pm$ should be replaced by ${\frak b}_\pm$ in the line preceding the
theorem.

{\bf Remark 2.} Another proof of Theorem 50 of \cite{EV1}
(different from the original one) can be obtained by
sending $\lambda$ to $\infty$ in the ABRR equation
(see \cite{ABRR,ES}).

Define $A_V':=A_V^+q^{\frac{h(h+2)}{4}}$. It follows from
Corollary \ref{limits}, part (iii) that on $V_m$,
one has $A_V'=A_V^\infty q^{m(m+2)/4}$.

\begin{corollary}\label{rmat} One has
$$
A_{V\otimes U}'={\mathcal R^{21}}(A_V'\otimes A_U').
$$
\end{corollary}

\begin{proof}
Straightforward from Proposition \ref{tensprod}.
\end{proof}

\begin{proposition}\label{EtoF} One has in $V$:
$$
A_V'f=-q^{-2}eA_V',\ A_V'e=-q^2fA_V'.
$$
and
$$
A_V^\infty f=-q^{-2}eA_V^\infty,\ A_V^\infty e=-q^2fA_V^\infty.
$$

\end{proposition}

\begin{proof} We prove the relations for $A_V'$; the relations
  for $A_V^\infty$ follow automatically since these two
operators are proportional.

If $V$ is 2-dimensional, then $A_V'$ is known, and
  it is straightforward to establish the result.
 From this and Corollary \ref{rmat} it follows
that the result is true if $V$ is the tensor product of any number
of 2-dimensional representations. But
any finite dimensional representation is contained in such a
product, so we are done.
\end{proof}

Now let $v_m$ be a highest weight vector of $V=V_m$.
Let us compute the operator of $A_V'$ in the basis
$v_{m-2j}:=\frac{f^j}{[j]_q!}v_m$, $j=0,...,m$.

\begin{proposition} One has
$$
A_V^\infty v_{m-2j}=(-1)^jq^{m-2j}v_{2j-m}.
$$
\end{proposition}

\begin{proof}
First of all observe that in $V_1^{\otimes m}$, one has
$\Delta_m(f^m)v_+^{\otimes m} =[m]_q!v_-^{\otimes m}$ (where
$\Delta_m$ is the iterated coproduct). On the other hand, using
the expression for $A_V^+$ for the 2-dimensional $V$, and the
expression for $A_{V_1\otimes V_2}^+$ given in Proposition
\ref{tensprod}, we get
$$
A_{V_1^{\otimes m}}^+(v_+^{\otimes m})=q^m(v_-^{\otimes m}).
$$
But $V_m$ is the submodule of
$V_1^{\otimes m}$ generated by
$v_+^{\otimes m}$.  Thus, for
$$
A_V^\infty v_m=A_V^+ v_m=q^mv_{-m}.
$$
Now the result follows from
Proposition \ref{EtoF} by a direct calculation.
\end{proof}

\subsection{The operators $B^\pm_V(\la)$}

Following \cite{TV},
define the operators $B^{\pm}_{V}(\la)$ by the formula
$$
B^{\pm}_{V}(\la)=(A_{V}^\pm)^{-1}A_{V}(\la).
$$

\begin{proposition} (i) $B^{\pm}_V$ preserve the weight
decomposition.

(ii) $B^+_V\to 1$ as $q^\la\to +\infty$;
$B^-_V\to 1$ as $q^\la\to 0$.

(iii) One has
$$
B_{V_m}^{\pm}(\la)|_{V_m[m-2k]}=
\prod_{j=1}^k\frac{(\la+1+j)_{q^{\mp 2}}}{(\la+1-m+2k-j)_{q^{\mp 2}}},
$$
where $(a)_q:=\frac{q^a-1}{q-1}$ is the nonsymmetric q-analog of
$a$.

(iv) If $q=1$, then $B^+_{V}$, $B^-_{V}$ are equal to
the operator $B_V$ from Corollary \ref{limits}.
\end{proposition}

\begin{proof} The proof of this proposition is straightforward
from the previous results.
\end{proof}

\begin{proposition} \label{Aforsl(2)}
The operator $B^+_V(\la)$ is given
 by the action in $V$ of the universal
element
$$
B^+(\la)=p(\la,h,e,f),
$$
where
$$
p(\la,h,e,f):=\sum_{k=0}^\infty q^{k(k-1)/2}\frac{q^{-k(\la+1)}}{[k]_q!}
f^ke^k\prod_{\nu=0}^{k-1}\frac{1}{[\la-h-\nu]_q}.
$$
In particular,
for $q=1$, the operator $B^\pm=B$ coincides with the operator
$B$ from section 2.5 of \cite{TV}.
\end{proposition}

{\bf Remark.} Similar formulas can be deduced for $B^-$ and $B$.

\begin{proof} This is proved by a straightforward calculation
  with intertwiners, which generalizes to the q-case
the calculation of \cite{TV}. Another proof is given as a remark
in Section 6.
\end{proof}

\section{The main construction
for any $\g$: the dynamical Weyl group}

\subsection{The operators $A_{w,V}(\la)$}

Let us return to the situation of a general $\g$.
Let $V$ be an integrable $U_q(\g)$-module,
and $\nu$ be a weight of $V$.
For any simple reflection
$s_i\in \W$, we define an operator-valued
rational function $A_{s_i,V}(\la): V[\nu]\to V[s_i\nu]$
by the formula
$$
A_{s_i,V}(\la)=A_{s,V'}(\la(h_i)),
$$
where $V'$ is the $U_{q_i}({{\frak sl}_2})$ submodule of $V$
generated by $V[\nu]$.

Let $w\in \W$, and let $w=s_{i_1}...s_{i_l}$ be a reduced
decomposition. Let us call this reduced decomposition $\delta$.

\vskip .05in

{\bf Definition:}
Define the operator $A_{w,V,\delta}(\la): V[\nu]\to
  V[w\nu]$
by
$$
A_{w,V,\delta}(\lambda)=A_{s_{i_1}}(s_{i_2}...s_{i_l}\cdot \la)...
A_{s_{i_{l-1}}}(s_{i_l}\cdot \la)A_{s_{i_l}}(\la).
$$
\vskip .05in

Thus,
the function $A_{w,V,\delta}(\lambda)|_{V[\nu]}$ uniquely extends
to a rational
operator valued function of the variables $(\lambda,\alpha_i)$, respectively
$q^{(\lambda,\alpha_i)}$.
This function is generically invertible.

The following two results play a crucial role
in our considerations.

Let $\Phi: M_\la\to M_\mu\otimes V$ be an intertwiner,
$\Phi v_\la=v_\mu\otimes <\Phi>+...$
Assume that $\la$ is dominant, and $\la(h_i)$ are sufficiently
large for all $i$.

\begin{proposition}\label{restr}  For any $\delta$ and $w\in \W$, one has
$$
\Phi v_{w\cdot \la,\delta}^\la=v_{w\cdot \mu,\delta}^\mu\otimes
A_{w,V,\delta}(\la)<\Phi>+\text{lower weight
  terms}.
$$
\end{proposition}

\begin{proof} The statement follows easily
from the definitions by induction on the length
of $w$.
\end{proof}

\begin{corollary}\label{indep}
The operator $A_{w,V,\delta}(\la)$ is independent
of the reduced decomposition $\delta$ of $w$.
\end{corollary}

\begin{proof} If $\la$ is large dominant,
the statement is clear
from Proposition \ref{restr}, since $v_{w\cdot \la,\delta}^\la$
is independent of
  the reduced decomposition $\delta$ of $w$ by the quantum Verma
identities (Lemma \ref{sing-v}).
For an arbitrary $\la$, the statement follows from the
  fact that a rational function is completely determined by its
  values at large dominant weights.
\end{proof}

Thus, we will denote $A_{w,V,\delta}(\la)$ simply by
$A_{w,V}(\la)$. Proposition
\ref{restr} shows that
as for $U_q({{\frak sl}_2})$, the operator $A_{w,V,\delta}(\la)$ is
the effect, at the level of expectation values,
of the restriction of an intertwiner from
$M_\la$ to $M_{w\cdot \la}$.

The following lemma is an  easy consequence of
the definition.

\begin{lemma}\label{Aww}
If $w_1,w_2\in\W$, $l(w_1w_2)=l(w_1)+l(w_2)$, then
\bea
A_{w_1w_2,V}(\la)\,=\,A_{w_1,V}(w_2\cdot\la)A_{w_2,V}(\la)\,
\eea
on $V[\nu]$.
\end{lemma}

More generally, we can define the operators
$A_{w,V}(\la)$ for any element $w$ of the braid group $\tilde \W$.
Namely, let $s_i^{-1}$ be the inverse of $s_i$ in
$\tilde\W$, and define the operator $A_{s_i^{-1},V}(\la)$ by
$$
A_{s_i^{-1},V}(\la)=A_{s_i,V}(s_i\cdot \la)^{-1}.
$$
Now, if $w=s_{i_1}^{\varepsilon_1}...s_{i_l}^{\varepsilon_l}$
where $\varepsilon_i=\pm 1$, then
we define $A_{w,V}(\la)$ by the formula
$$
A_{w,V}(\lambda)=A_{s_{i_1}^{\varepsilon_1}}(s_{i_2}...s_{i_l}\cdot \la)...
A_{s_{i_{l-1}}^{\varepsilon_{l-1}}}(s_{i_l}\cdot \la)
A_{s_{i_l}^{\varepsilon_l}}(\la).
$$
It is easy to see from
Corollary \ref{indep} or Lemma \ref{Aww} that this definition
is independent on the product representation of $w$ but depends
only of $w$ itself,
so the notation $A_{w,V}$ is non-ambiguous.

\subsection{The dynamical Weyl group}

The above results imply that using the operators
$A_{w,V}(\lambda)$,
one can define a $\C$-linear action of the braid group $\tilde\W$
on the space of meromorphic functions of $\lambda$ with values in
$V$, by the formula
$$
(w\circ f)(\lambda)=A_{w,V}(w^{-1}\cdot \lambda)f(w^{-1}\cdot
\lambda), w\in \tilde \W
$$

Now we will give the main definitions of this paper.

\vskip .05in
{\bf Main definition 1:} Let $V$ be an integrable
$U_q(\g)$-module. The
$\tilde\W$-action $f\to w\circ f$ on $V$-valued functions on
$\h^*$ will be called
{\bf the shifted
dynamical action}.
\vskip .05in

We would also like to give a name to the operators $A_{w,V}(\la)$,
since they play a central role in the paper.
We defined these operators for all $w\in \tilde\W$. However,
the operators $A_{w,V}$ corresponding to elements $w$ 
of $\W$ (regarded as a subset of $\tilde \W$)
are especially remarkable and occur especially often in applications.
Therefore, we will restrict our attention to them, and make the
following definition.

\vskip .05in
{\bf Main definition 2:} We call
the collection of operator valued rational functions
$A_{w,V}(\lambda)$, $w\in \W$
{\bf the dynamical Weyl group} of $V$.
\vskip .05in

{\bf Remark.} One of the important properties of dynamical
R-matrices is that they tend to usual R-matrices when the
dynamical parameter $\lambda$ goes to infinity
(see \cite{EV1}, Theorem 50).
On the other hand, we will show later
that when $\lambda$ goes to infinity, the operators
$A_{w,V}(\lambda)$ tend to elements of the usual classical or
quantum Weyl group acting in $V$. This justifies the term
``dynamical Weyl group''.

We also define the (unshifted) {\bf dynamical action} $w\to w*$
of $\tilde \W$ on functions of $\lambda$ with values in $V$ as
follows. We introduce the operators
$$
\A_{w,V}(\la)=A_{w,V}(-\la-\rho+\frac{1}{2}h).
$$
Then by the definition
$$
(w*f)(\la)=\A_{w,V}(w^{-1}\la)f(w^{-1}\la).
$$

{\bf Remark.} Observe that for $w_1,w_2\in\W$, such that $l(w_1w_2)=
l(w_1)+l(w_2)$, one has
$$
\A_{w_1w_2,V}(\la)=\A_{w_1,V}(w_2\la)\A_{w_2,V}(\la).
$$
This implies that $w*$ is indeed an action of $\tilde \W$.

\subsection{The dynamical Weyl group of
the tensor product and the dual representation}\label{tens}

The following is one of the main properties of the dynamical Weyl
group.

\begin{lemma}\label{com-A}
Let $U, V$ be integrable $U_q(\g)$-modules. Let $w\in \W$. Then
\bean\label{Com-A}
A_{w,U\T V}(\la)J_{UV}(\la)\,=\,J_{UV}(w\cdot\la)
A_{w,V}^{(2)}(\la)A_{w,U}^{(1)}(\la-h^{(2)})
\,
\eean
where $A^{(1)}$ denotes $A\otimes 1$,
and $A^{(2)}$ denotes $1\otimes A$.
\end{lemma}

\begin{proof} As for $U_q({{\frak sl}_2})$, the lemma is an easy consequence of the
  definitions: it expresses the fact that the operation of
fusion of intertwiners
  commutes with the operation of restriction of intertwiners to submodules.
\end{proof}

\begin{corollary}\label{Rtrans}
For any $w\in \W$ the dynamical R-matrix $R_{VU}(\lambda)$
  satisfies the relation
$$
R_{VU}(w\cdot \lambda)=
A^{(2)}_{w,U}(\lambda)
A^{(1)}_{w,V}(\lambda-h^{(2)})R_{VU}(\lambda)
A^{(2)}_{w,U}(\lambda-h^{(1)})^{-1}
A^{(1)}_{w,V}(\lambda)^{-1}
$$
\end{corollary}

{\bf Remark.}
Here $A^{(2)}_{w,U}(\lambda-h^{(1)})^{-1}$
is the inverse of the operator $A^{(2)}_{w,U}(\lambda-h^{(1)})$.

Now let $\g$ be finite dimensional, and let us calculate
the dynamical Weyl group operators on the dual to a given representation.

Recall that the dual space $U^*$ of a $U_q(\g)$-module $U$
has two module structures. The first one is given by  $a\to (S(a)|_{U})^*$
and the second one is given by  $a\to (S^{-1}(a)|_{U})^*$. The corresponding
$U_q(\frak g)$-modules are denoted $U^*$ and ${}^*U$. These modules are
isomorphic, via $q^{2\rho}:{}^*U\to U^*$.

Let $Q_U(\lambda): U\to U$ be given by
$Q_U(\la)=\sum (c_i')^*c_i$, where $\sum c_i\otimes c_i'=
J_{U,{}^*U}(\lambda)$ (this operator was introduced in \cite{EV1,EV2}).

\begin{proposition}\label{Adual} For any $w\in \W$ one has
$$
A_{w,{}^*U}(\lambda)^*=Q_U(\lambda)A_{w,U}(\lambda-h)^{-1}Q_U(w\cdot \la)^{-1},
$$
and hence
$$
A_{w,U^*}(\lambda)^*=q^{2\rho}Q_U(\lambda)A_{w,U}(\lambda-h)^{-1}
Q_U(w\cdot \la)^{-1}q^{-2\rho}.
$$
\end{proposition}

{\bf Remark.} This formula has recently been used in the theory
of trace functionals on non-commutative moduli spaces of flat
connections,
see \cite{RS}.

\begin{proof} It is enough to establish the first formula.
The proof of this formula
is obtained by specializing formula (\ref{Com-A})
to the case $V={}^*U$,
dualization of the second component, and
multiplication of the components.
\end{proof}

Here is another formula for
the dynamical Weyl group of the dual representation,
which is valid on the zero weight subspace,
and involves a sign change for $\la$.

\begin{proposition}\label{dualonweight0} One has
$$
\A_{w,U^*}(w^{-1}\la)^*|_{U[0]}=\A_{w,{}^*U}(w^{-1}\la)^*|_{U[0]}=
\A_{w^{-1}}(-\la)|_{U[0]}
$$
\end{proposition}

\begin{proof} Let $w=s_{i_1}...s_{i_l}$ be a reduced
decomposition. Using the Main Definition,
we have
$$
\A_{w,U^*}(w^{-1}\la)=
\A_{s_{i_1},U^*}(s_{i_2}...s_{i_l}w^{-1}\la)...
\A_{s_{i_l},U^*}(w^{-1}\la)=
\A_{s_{i_1},U^*}(s_{i_1}\la)...
\A_{s_{i_l},U^*}(w^{-1}\la)=
$$
$$
\A_{s_{i_1},U^*}(-\la)...
\A_{s_{i_l},U^*}(-s_{i_{l-1}}...s_{i_1}\la)
$$
(in the last equality we used that $A_{s_i,V}(\la)$ depends only
of $(\la,\alpha_i)$). This implies that
$$
\A_{w,U^*}(w^{-1}\la)^*=
\A_{s_{i_l},U^*}(-s_{i_{l-1}}...s_{i_1}\la)^*
...\A_{s_{i_1},U^*}(-\la)^*
$$
But for $\g={{\frak sl}_2}$ we obviously have
$$
\A_{s,U^*}(\la)^*=\A_{s,U}(\la)
$$
on $U[0]$. This, together with the Main Definition
implies the result.
\end{proof}

\section{Further properties of the dynamical Weyl group}

\subsection{Limits of $A_{w,V}(\la)$ for any $\g$}

Let $\g$ be any Kac-Moody algebra, and $w\in \W$.
Let $V$ be an integrable $U_q(\g)$-module.
For 
a root $\alpha$, let $\varepsilon_{\alpha}(C)$ be the sign of 
$(\la,\alpha)$ in a Weyl chamber $C$. 

\begin{proposition} \label{limits1}
(i) If $q=1$, then the map $A_{w,V}(\lambda)$
has a limit $A_{w,V}^\infty$ at $\lambda\to\infty$
in a generic direction.
If $q\ne 1$, then 
for any Weyl chamber $C$ the map $A_{w,V}(\lambda)$ has
a limit $A_{w,V}^C$ as
$\la\to \infty$ in a generic direction in the chamber $C$.

(ii) Let us agree that for $q=1$, $A_{w,V}^C:=A_{w,V}^\infty$
for all $C$. Then, for 
every reduced decomposition $w=s_{i_1}...s_{i_l}$,  
and any $q$, one has
$$
A_{w,V}^C=A_{s_{i_1},V}^{\varepsilon_{\alpha^1}(C)}...A_{s_{i_l},V}^{
\varepsilon_{\alpha^l}(C)},
$$
where $A_{s_i,V}^\pm $ are the elements $A_{s,V}^\pm$ for the
subalgebras generated by $e_i,f_i,h_i$ (for $q=1$) or
$e_i,f_i,q^{bh_i}$ (for $q\ne 1$). 
\end{proposition}

\begin{proof}
The proof follows easily from the results of the
previous section and the Main Definition. 
\end{proof}

We will mostly consider the operators $A_{w,V}^C$ if 
$C$ is the dominant or the antidominant chamber. 
In this case we will denote $A_{w,V}^C$ by $A_{w,V}^+$ and $A_{w,V}^-$, 
respectively.

Part (ii) of Proposition \ref{limits1} implies the following
statements

\begin{corollary}\label{limits1.1} The assignments
$s_i\to A_{s_i,V}^+$, $s_i\to A_{s_i,V}^-$ extend to actions of
$\tilde\W$ on $V$.
\end{corollary}

\begin{proof} Clear. \end{proof}

\begin{corollary}\label{limits1.2}
On $V[0]$, one has $A_{w,V}^-A_{w^{-1},V}^+=1$ for any
$w\in \W$.
\end{corollary}

\begin{proof} To prove the statement for any
Lie algebra, it is enough to do so for $\g={{\frak sl}_2}$, in
which case the result is an easy consequence of the results of
Section 3.
\end{proof} 

\subsection{The relation to the quantum Weyl group}

In this subsection we will discuss the relationship of the dynamical
Weyl group with the quantum Weyl group for $U_q(\g)$.

The quantum Weyl group was introduced by Soibelman
for finite dimensional $\g$ (see \cite{KoSo})
and by Lusztig in the Kac-Moody case
(see \cite{Lu}).
It is discussed in many books and papers,
which use various (though pairwise equivalent) conventions.
We will use the conventions of the paper \cite{Sa}.

The quantum Weyl group element $\Bbb S$ in a completion of
$U_q({{\frak sl}_2})$ is defined by the formula \cite{Sa}:
$$
\Bbb S=\exp_{q^{-1}}(q^{-1}eq^{-h})\exp_{q^{-1}}
(-f)\exp_{q^{-1}}(qeq^h)q^{h(h+1)/2},
$$
where $\exp_p(x)=\sum_{m\ge 0}\frac{p^{m(m-1)/2}}{[m]_p!}x^m$.
According to \cite{Sa}, Proposition 1.2.1, this element
acts in $V=V_m$ as
$$
\Bbb S v_{m-2j}=(-1)^{m-j}q^{(m-j)(j+1)}v_{2j-m}
$$
Therefore, in $V$ we have:
$$
\Bbb S|_{V[m-2j]}=(-1)^m q^{(m-j+1)j}A_V^\infty=(-1)^mA_V^+.
 $$

\begin{proposition}\label{AS}
For a finite dimensional $U_q({{\frak sl}_2})$-module $V$, one has
$$
A_V^+=(-1)^h\Bbb S ,\quad A_V^-=q^h\Bbb S^{-1}
$$
 where $(-1)^h$ is the
transformation acting by $-1$ on even-dimensional
irreducible modules, and by
$1$ on odd dimensional ones
(i.e. the quantum analog of the element $-1$ of the group
$SL(2)$).
\end{proposition}

\begin{proof} This follows at once
by comparing the actions of both sides on basis vectors of $V$.
\end{proof}

Now let $\g$ be a Kac-Moody algebra.
Following \cite{Sa}, define the element
$S_i$ to be the element $\Bbb S$ of the
simple root subalgebra of $U_q(\g)$
corresponding to the simple root $\alpha_i$.
The elements $S_i$ define operators on any integrable
module (the ``symmetries of an integrable module'' defined
by Lusztig \cite{Lu}).

As an immediate corollary of Proposition \ref{AS}, we obtain
the following well known and important result
(see \cite{Lu},\cite{Sa}).

\begin{corollary} The elements $S_i$
satisfy the braid relations of $\tilde \W$ in any integrable module.
\end{corollary}

\begin{proof} The result
follows immediately from
the equality $A_V^-=q^h\Bbb S^{-1}$
and Proposition \ref{limits1}.
\end{proof}

{\bf Remark.} In particular, this result is valid
for $q=1$ and yields a braid group action on
integrable modules over a classical Kac-Moody algebra.
This action factorizes through an extension of $\W$
by a group isomorphic to $(\Bbb Z/2)^r$.

For $w\in s_{i_1}^{\varepsilon 1}...
s_{i_l}^{\varepsilon_l}\in \tilde \W$, let $\Bbb S_w$ be the quantum
Weyl group operator corresponding to $w$ (i.e.
$\Bbb S_w=S_{i_1}^{\varepsilon_1}...S_{i_l}^{\varepsilon_l}$.

\begin{proposition} For an integrable
$U_q(\g)$-module $V$ and $w\in \tilde \W$,
one has $A_{w,V}^+=(-1)^{\rho^\vee
    -w\rho^\vee}\Bbb S_w$,
where $\rho^\vee$ is an element of $\h$
such that $\alpha_i(\rho^\vee)=1$ for all
$i$.
\end{proposition}

\begin{proof} This is immediate from Proposition \ref{AS} and the
  definition.
\end{proof}

{\bf Remark.} 
Analogously to this proposition, 
the limits $A_{w,V}^C$ of the operators
$A_{w,V}(\la)$ in Weyl chambers $C$ other than the dominant one 
give rise to
other variants of the quantum Weyl group.

\subsection{The dynamical action of the pure braid group}
Define the {\bf pure
  braid group} $P\tilde\W$ to be
the kernel of the natural homomorphism $\tilde\W\to \W$. This
group is the smallest normal subgroup of $\tilde \W$ which
contains $s_i^2$ for all $i$. So it is generated by the elements
$w s_i^2 w^{-1}, w\in \tilde \W$.

The dynamical action of $\tilde \W$ on functions of $\la$
with values in an integrable
$U_q(\g)$-module $V$ induces an action of $P\tilde \W$
on functions with values in any weight subspace $V[\nu]$ of $V$,
which is linear over the field $F$ of scalar meromorphic
functions of $\lambda$.

\begin{proposition}\label{scalar} The group $P\tilde \W$ acts in
$V[\nu]$ by multiplication by a character $\chi: P\tilde \W\to F^*$,
which is $\tilde \W$-equivariant
(i.e. $\chi(wpw^{-1})(\la)=\chi(p)(w\la)$),
and is determined by
the relations
$$
\chi(s_i^2)=(-1)^{\nu(h_i)}\frac{[(\la+\nu,\alpha_i)]_q}{[\la]_q},
$$
\end{proposition}

\begin{proof} It is sufficient to prove the result for ${{\frak sl}_2}$.
In this case, we have a unique simple reflection $s$, and
$$
(s*)^2=\A_V(-\la)\A_V(\la)=(A_V^\infty)^2
B_V(\la-1)B_V(-\la-1)_{V[\beta]}.
$$
But $(A_V^\infty)^2=(-1)^h$ on $V$. Therefore, by the explicit formula
for $B_V$ we have $(s*)^2=(-1)^{\beta(h)}\frac{[\la+\beta]_q}{[\la]_q}$,
as claimed.
\end{proof}

It is easy to see that if $V$ is an integrable module over
$U_q(\g)$ then the dynamical Weyl group operators preserve the
space of functions on $\h^*$ with values in the zero weight
subspace $V[0]$. So let us consider the dynamical action of the
braid group restricted to this space.

\begin{corollary}\label{squareisone} The dynamical action of
the pure braid group $P\tilde \W$ on
functions of $\la$ with values in $V[0]$
is trivial. Therefore,
the dynamical action of $\tilde \W$ on this space
induces an action of the Weyl group  $\W$.
In particular, the 1-cocycle relation
$A_{w_1w_2,V}(\la)=A_{w_1,V}(\la)A_{w_2,V}(w_1\cdot \la)$
of Lemma \ref{Aww} is satisfied for any $w_1,w_2\in \W$
(i.e. without the requirement $l(w_1w_2)=l(w_1)+l(w_2)$).
\end{corollary}

\begin{proof} Follows immediately from Proposition
\ref{scalar}.
\end{proof}

{\bf Remark.} Note that these results are false
for the usual (non-dynamical)
quantum Weyl group, unless $q=1$. The action of the pure part of
the quantum Weyl group does not, in general, reduce to
a character, and the action on the zero weight subspace is, in
general, nontrivial. This seemingly paradoxical situation
is similar to the situation with trigonometric R-matrices
with a spectral parameter: these R-matrices satisfy the involutivity
condition and hence define representations of $S_n$, but tend at
infinity to R-matrices without spectral parameter, which fail
to satisfy involutivity and hence define only braid group
representations.

\subsection{The operators $B^{\pm}_{w,V}(\la)$}

Let $w\in \W$.
Following \cite{TV},
define the operators $B^{\pm}_{w,V}(\la)$ by the formula
$$
B^{\pm}_{w,V}(\la)=(A_{w,V}^\pm)^{-1}A_{w,V}(\la).
$$
For example, in the $U_q({{\frak sl}_2})$ case,
and $w$ being the only nontrivial element,
we have $B_{w,V}^\pm(\la)=B_V^\pm(\la)$
(see the notation in Section 3).

\begin{proposition} (i) $B^{\pm}_{w,V}$ preserves the weight
decomposition.

(ii) We have $B_{w,V}^+\to 1$ as $q^{(\la,\alpha_i)}\to +\infty$;
$B_{w,V}^-\to 1$ as $q^{(\la,\alpha_i)}\to 0$.

(iii) If $q=1$, then $B^+_{w,V}$ equals $B^-_{w,V}$.
\end{proposition}

\begin{proof} The proof of this proposition is straightforward.
\end{proof}

Thus, if $q=1$, we will denote $B_{w,V}^+=B_{w,V}^-$ simply by
$B_{w,V}$.

From the definition of the dynamical Weyl group it follows that
the operators $B_{w,V}^\pm$, like the operators $A_{w,V}$, admit
a factorization into a noncommutative product of $l(w)$ terms.
For the sake of brevity, we will consider only the operators
$B_{w,V}^+$; the story for $B_{w,V}^-$ is analogous.

Let $w\in \W$. Fix a reduced decomposition $\delta$ of $w$:
$w=s_{i_1}...s_{i_l}$.
Recall that in Section 4 we assigned to this reduced decomposition
a sequence of roots $\alpha^1,...,\alpha^l$.

It follows from Lusztig's theory of braid group actions on
$U_q(\g)$ \cite{Lu} that there exists a unique element
$e^\delta_{\pm \alpha^j}$ of $U_q(\g)$ which acts in any
integrable module $V$ as
$$
e^\delta_{\pm \alpha_j}=
A_{s_{i_l}...s_{i_{j+1}}}^+e_{\pm
\alpha_{i_j}} (A_{s_{i_l}...s_{i_{j+1}}}^+)^{-1}
$$
(a ``Cartan-Weyl
generator''). We also let
$h_{w\alpha_i}=w(h_i)$.

\begin{proposition}\label{produ}
The operator $B^+_{w,V}(\la)$ is obtained
 by the action in $V$ of the universal
element
$$
B^+_w(\la)=
\prod_{j=1}^l p((\la+\rho)(h_{\alpha^j})-1, h_{\alpha^j},e_{\alpha^j}^\delta,
e_{-\alpha^j}^\delta),
$$
where the indices increase from left to right.
\end{proposition}

\begin{proof} Straightforward from the definition.
\end{proof}

In particular, Proposition \ref{produ} implies that the operator
$B^+_{w,V}(\la)$, unlike the operator $A_{w,V}(\la)$, is well
defined on weight spaces of any module $V$ over $U_q(\g)$ from
category $\mathcal O$ (not necessarily integrable). Namely, it is
defined by the product formula of Proposition \ref{produ}.

Proposition \ref{produ} also implies the following determinant
formula for $B^+_{w,V}$ acting on a weight subspace.

Let $w\in \W$ and $w=s_{i_1}...s_{i_l}$ be a reduced decomposition of $w$.
Let $\beta$ be a weight
of a finite dimensional $U_q(\g)$-module $V$.
Let
$$
B^+_{\alpha\beta k}(q,\la)=
\prod_{j=1}^{k}\frac{q^{-2((\la+\rho,\alpha^\vee)+j)}-1}
{q^{-2((\la+\rho-\beta,\alpha^\vee)-j)}-1}.
$$

\begin{corollary}\label{deter} One has
$$
\text{det}(B^+_{w,V}(\la)|_{V[\beta]})=
\prod_{j=1}^l
\prod_{k\ge 0}B^+_{\alpha^j\beta k}(q_{i_j},\la)^{d_V({\alpha^j},\beta,k)},
$$
where
$$
d_V(\alpha,\beta,k)=\dim(V[\beta+k\alpha])-
\dim(V[\beta+(k+1)\alpha]).
$$
\end{corollary}

\begin{proof} Straightforward from the definition.
\end{proof}

\subsection{The operators $B^+_{w_0,V}(\la)$ and extremal projectors}

Let $\g$ be finite dimensional, $w_0$ the maximal element of the
Weyl group of $\g$. Let $M_\mu$ be the Verma module with highest
weight $\mu$. Let $\gamma$ be a nonnegative linear combination of
simple roots.

\begin{proposition} \label{extproj}
For generic $\mu$, the operator-valued function
$B^+_{w_0,M_\mu}(\la)|_{M_\mu[\mu-\gamma]}$ is regular at the
point  $\la=-2\rho$, and the operator $B^+_{w,M_\mu}(-2\rho)$ is
the projector to the highest weight space of $M_\mu$ along the
sum of other weight spaces (the extremal projector).
In other words, $B^+_{w,M_\mu}(-2\rho)|_{M_\mu[\mu-\gamma]}=\delta_{\gamma 0}
\text{Id}$.
\end{proposition}

\begin{proof} The first statement is immediate
from Proposition \ref{produ}, as for
$\g={\frak sl}_2$, the operator valued function $B^+(\la)$
of $\la$
is regular at integer values of $\la$ when restricted to
a weight subspace of a generic Verma module.

Let us now prove the second statement.
Because $\mu$ is generic, any
nonzero homogeneous vector $v\in M_\mu$
of weight different from $\mu$ is not singular.
Thus, it suffices to show that for any non-singular homogeneous
vector $v$, one has $B^+_{w_0,M_\mu}v=0$.

Let $v$ be a non-singular homogeneous vector in $M_\mu$. Then
there exists an index $i$ such that $e_iv\ne 0$. We will assume
that $v$ is an eigenvector for the Casimir operator of
$U_{q_i}({\frak sl}_2)$, since this does not cause a loss of
generality. Let $Y_v$ be the submodule of $M_\mu$ over
$U_{q_i}({\frak
  sl}_2)$ generated by $v$. Then $Y_v$ is a Verma module,
and $v$ is a nonzero homogeneous vector in $Y_v$ which is not
singular.

Let $w_0=s_{i_1}...s_{i_l}$ be the reduced decomposition of
$w_0$, such that $i_l=i$. By Proposition \ref{produ},
to this decomposition there corresponds a factorization
$$
B_{w_0,M_\mu}^+(-2\rho)=\Pi\cdot p(-2,h_i,e_i,f_i),
$$
where $\Pi$ is the product of the terms in the product formula
for $B^+_{w_0,V}$ corresponding to all but the last factor in the
reduced decomposition of $w_0$. Thus, it suffices to show that
$p(-2,h_i,e_i,f_i)v=0$.

But this is immediate from the product formula for the operator
$B$ in Corollary \ref{limits}: the first factor in this product
has numerator $\la+2$, and so for $\la=-2$, the product vanishes
whenever the set of indices over which the product is taken is
nonempty. The proposition is proved.
\end{proof}

We note that extremal projectors appeared in \cite{AST}
and there is an extensive theory of them (see
e.g. \cite{Zh1,Zh2}). In particular, Proposition \ref{extproj},
in the case of finite dimensional Lie algebras and $q=1$,
essentially appears in \cite{Zh1,Zh2} (see Theorem 2 of \cite{Zh2}).

\subsection{The dynamical Weyl group of a locally finite module}

Let $V$ be an $\h$-diagonalizable $U_q(\g)$-module.
Recall that $V$ is said to be locally finite
if any vector of $v\in V$ generates a finite dimensional
submodule over the quantum $U_{q_i}({{\frak sl}_2})$
subalgebras corresponding to all
simple roots.

\begin{lemma}\label{fin} Let
$\g$ be finite dimensional, and
let $V$ be a locally finite $U_q(\g)$-module. Then
the submodule $Y_v$ generated in $V$
by any vector $v$ is finite dimensional.
\end{lemma}

\begin{proof} It suffices to assume that $v$ is homogeneous
with respect to the weight decomposition.

Let $\lbrace{ e_\alpha\rbrace} $ be the Cartan-Weyl generators of
$U_q(\g)$ corresponding to some reduced decomposition of the
maximal element of $\W$. By the PBW theorem (see \cite{Sa}), the
submodule $Y_v\subset V$ is given by
$Y_v=\prod_{\alpha}\C[e_\alpha] v$, where $\C[e_\alpha]$ is the
algebra of polynomials of $e_\alpha$, and the product is taken
over all roots in a suitable order. Since the product is finite,
and $V$ is a sum of finite dimensional $\C[e_\alpha]$-modules
(because $e_\alpha$ is conjugate to some $e_i$ under Lusztig's
braid group action on $U_q(\g)$), we have that $Y_v$ is finite
dimensional.
\end{proof}

It is clear that the operators $A_{w,V}(\la)$ are
  well defined for any locally finite module $V$
over $U_q(\g)$.

Therefore, for any such $V$ we can define
the operators $w\circ$, $w\in \tilde\W$,
on $V$-valued functions of $\la$ as before.
The above lemma allows us to prove the following.

\begin{corollary} For any locally finite $U_q(\g)$-module $V$, the
operators $w\circ $ define an action of $\tilde \W$.
\end{corollary}

\begin{proof} We have to check that the operators $s_i\circ$
  satisfy the braid
  relations. This can be checked on finite dimensional
Lie algebras of rank 2. But in this case, according
to Lemma \ref{fin}, everything reduces to the case
when $V$ is finite dimensional, where the statement is known.
\end{proof}

\section{The dynamical Weyl group element
corresponding to the maximal element $w_0$ of $\Bbb W$}

\subsection{The expression for $A_{w_0,V}$ and $B_{w_0,V}^+$}

In this section we will study the operator $A_{w_0,V}$ for the
maximal element $w_0$ of $\W$ in the case of finite dimensional
Lie algebras. This operator is especially important because, as
we will show below, it is closely related to the operator
$Q_V(\la)$, which is (as we will show elsewhere) the matrix of
inner products of trace functions, and therefore generalizes
squared norms of Macdonald polynomials. The material of this
section allows one to calculate explicitly the operator $Q_V$ and
its determinant (see below), and therefore prove a matrix analog of the
Macdonald inner product identities.

So, let $\g$ be finite dimensional, and let $w_0\in \Bbb W$
be the maximal element.
Let $\{x_i\}$ be an orthonormal basis of $\h$.
Let $U,V$ be finite dimensional modules over $U_q(\g)$.

\begin{lemma}\label{oppo}
One has in $V\otimes U$:
$$
A_{w_0,V\otimes U}^+={\mathcal R}_0^{21}(A_{w_0,V}^+\otimes A_{w_0,U}^+),
$$
and
$$
\text{Ad}(A_{w_0,V}^+\otimes A_{w_0,U}^+)(\mathcal R)=
q^{\sum x_i\otimes x_i}\mathcal
R_0^{21}=q^{\sum x_i\otimes x_i}\mathcal
R^{21}q^{-\sum x_i\otimes x_i}.
$$
\end{lemma}

\begin{proof} The
statements are obtained from Lemma \ref{com-A} and
Corollary \ref{Rtrans} by passing sending $\la$ to
infinity. Namely, the
first equation follows from Lemma \ref{com-A}
and Theorem 50 of \cite{EV1}.
The second equation follows from
Corollary \ref{Rtrans} and Theorem 50 of \cite{EV1}.
\end{proof}

{\bf Remark.} Since $A_{w,V}^+$ are, essentially, the quantum
Weyl group operators, Lemma \ref{oppo} can also be deduced from
the theory of the quantum Weyl group
(see \cite{KoSo}). In fact, this lemma is nothing but the
coproduct property of the maximal element, which is an important
property of the quantum Weyl group. This property is instrumental in
deriving the Levendorski-Soibelman product formula for the
universal R-matrix (\cite{KoSo}).

Let $Q^\dagger_V(\la)=Q_{V^*}(\la)^*$
(in the notation of \cite{EV2},
$Q^\dagger_V=S(Q)|_V$, where $S$ is the antipode).

In this subsection we will prove the following
result.

\begin{thm}\label{Qdagger} One has
$$
A_{w_0,V}(\la)=A_{w_0,V}^+ Q^\dagger_V(\la).
$$
In other words, $Q^\dagger_V(\la)=B^+_{w_0,V}(\la)$.
\end{thm}

{\bf Remark.} Note that in the formula for $Q^\dagger$
obtained from the definition, the f-terms come on the left from the e-terms,
while in the product formula for $B_{w_0}^+$, the f-terms and the
e-terms are mixed with each other. Thus, the theorem provides a
way to perform a ``normal ordering'' of the terms in the product
formula for $B_{w_0}^+$.

The proof of Theorem \ref{Qdagger} occupies the rest of the
subsection.

Consider the universal fusion matrix $J(\la)$, i.e.
the element of the completed tensor square of $U_q(\g)$ which
acts in the product $V\otimes U$ of finite dimensional
$U_q(\g)$-modules as $J_{VU}(\la)$ (see \cite{EV1}).
Let $Q(\la)=\sum S^{-1}(c_i') c_i$, where $J(\la)=\sum c_i\otimes c_i'$,
and $S$ is the antipode. Let $Q^\dagger(\la)=S(Q(\la))$.
Then $Q|_V=Q_V$, $Q^\dagger|_V=Q^\dagger_V$.

It follows from formula (2.38) in \cite{EV2} that
$$
\Delta(Q(\la))=(S\otimes S)(J^{21}(\la))^{-1}
\cdot (Q(\la)\otimes Q(\la+h^{(1)}))
\cdot J(\la+h^{(1)}+h^{(2)})^{-1}.
$$
Applying the antipode to this equation and permuting the
components, we obtain
$$
\Delta(Q^\dagger(\la))=(S\otimes S)(J^{21}(\la+h^{(1)}+h^{(2)}))^{-1}
\cdot (Q^\dagger(\la-h^{(2)})\otimes Q^\dagger(\la))\cdot J(\la)^{-1}.
$$
Here we use the fact that $S^2=\text{Ad}(q^{2\rho})$ and the zero
weight property of $J$; note also that we replaced $h^{(i)}$ by
$-h^{(i)}$ in the $Q$-terms of the equation, since the antipode
changes the sign.

On the other hand, we have shown in Lemma \ref{oppo} that
$$
A_{w_0,V\otimes U}^+={\mathcal R_0^{21}}(A_{w_0,V}^+\otimes A_{w_0,U}^+).
$$
Therefore, we have the following lemma.

\begin{lemma}\label{AQ}
$$
A_{w_0,V\otimes U}^+
Q^\dagger_{V\otimes U}(\la)=
$$
$$
{\mathcal R_0^{21}}\cdot (A_{w_0,V}^+\otimes A_{w_0,U}^+)\cdot
(S\otimes S)(J^{21}(\la+h^{(1)}+h^{(2)}))^{-1}
\cdot (Q^\dagger(\la-h^{(2)})\otimes Q^\dagger(\la))
\cdot J(\la)^{-1}.
$$
\end{lemma}

We will also need the following lemma.

\begin{lemma}\label{twistsequal}  One has on $V\otimes U$:
$$
{\mathcal R_0^{21}}\cdot \text{Ad} (A_{w_0,V}^+\otimes A_{w_0,U}^+)
\cdot \left( (S\otimes S)(J^{21}(\la+h^{(1)}+h^{(2)}))^{-1}\right)=
J(w_0\cdot \la).
$$
\end{lemma}

\begin{proof}
We will first transform the equality to a convenient form, and
then show that both sides satisfy the same ABRR equation
(\cite{ABRR,ES}), which
has a unique solution. This will imply that the two sides are
equal.

Let us make a change of variable
$\Bbb J(\la)=J(-\la-\rho)$. Then the equation to be proved
takes the form
$$
{\mathcal R_0^{21}}\cdot\text{Ad}(A_{w_0,V}^+\otimes A_{w_0,U}^+)\cdot
\left( (S\otimes S)(\Bbb J^{21}(\la-h^{(1)}-h^{(2)}))^{-1}\right)=
\Bbb J(w_0\la).
$$
Let ${\mathcal J}(\la)=\Bbb
J(\la-\frac{1}{2}(h^{(1)}+h^{(2)}))$.
Then the equation takes the form
$$
{\mathcal R_0^{21}}\cdot\text{Ad}(A_{w_0,V}^+\otimes A_{w_0,U}^+)
\cdot\left( (S\otimes S)({\mathcal J}^{21})^{-1}\right)
(\la+\frac{1}{2}w_0
(h^{(1)}+h^{(2)}))=
{\mathcal J}(w_0\la+\frac{1}{2}(h^{(1)}+h^{(2)})).
$$
Replacing $\la+\frac{1}{2}w_0(h^{(1)}+h^{(2)})$
with $\la$, we obtain the equation
$$
{\mathcal R_0^{21}}\cdot\text{Ad}(A_{w_0,V}^+\otimes A_{w_0,U}^+)
\cdot\left( (S\otimes S)({\mathcal J}^{21}(\la))^{-1}\right)=
{\mathcal J}(w_0\la).
$$
To establish this equation, let us recall that
by Lemma 2.4 of \cite{EV2} (see also \cite{ABRR}),
the element $X(\la)={\mathcal
  J}(\la)$ satisfies the ABRR equation
$$
{\mathcal R^{21}}(q^{2\la})_1X(\la)=
X(\la)q^{\sum x_i\otimes x_i}(q^{2\la})_1.
$$
Therefore, the element
$Y(\la)=(S\otimes S)({\mathcal J}^{21}(\la))^{-1}$
satisfies the equation
$$
{\mathcal R}^{-1}(q^{2\la})_2Y(\la)=Y(\la)q^{-\sum x_i\otimes
  x_i}(q^{2\la})_2.
$$
Thus, using Lemma \ref{oppo},
we get that the operator on $V\otimes U$ given by
\linebreak $Z(\la)=
\text{Ad}(A_{w_0,V}^+\otimes A_{w_0,U}^+)(Y(\la))$
satisfies the equation
$$
({\mathcal R_0^{21}})^{-1}
q^{-\sum x_i\otimes x_i}(q^{2w_0\la})_2Z(\la)=Z(\la)q^{-\sum x_i\otimes
  x_i}(q^{2w_0\la})_2.
$$
Therefore, the operator $T(\la)={\mathcal R}_0^{21}Z(\la)$
satisfies
$$
(q^{2w_0\la})_2{\mathcal R^{21}}^{-1}
T(\la)=T(\la)q^{-\sum x_i\otimes
  x_i}(q^{2w_0\la})_2.
$$
Transforming this using the weight zero property of $T$, we get
$$
{\mathcal R}^{21}(q^{2w_0\la})_1T(\la)=
T(\la)(q^{2w_0\la})_1q^{\sum x_i\otimes x_i}.
$$
Now we note that the same equation is satisfied by ${\mathcal
J}(w_0\la)$, by Lemma 2.4 of \cite{EV2}. Both of these solutions
are triangular, with the diagonal part equal to $1$. But Lemma 2.4
of \cite{EV1} claims that such a solution is unique.

Therefore,
$T(\la)={\mathcal J}(\la)$, and the lemma is proved.
\end{proof}

\begin{corollary}
The operators
$E_V(\la):=A_{w_0,V}(\la)^{-1}A_{w_0,V}^+Q_V^\dagger(\la)$ have
the ``grouplike'' property
$$
E_{V\otimes U}(\la)=J(\la)(E_V(\la-h^{(2)})\otimes
E_U(\la))J^{-1}(\la).
$$
\end{corollary}

\begin{proof} This follows directly from
Lemma \ref{com-A}, Lemma \ref{AQ}, and Lemma \ref{twistsequal}.
\end{proof}

Now we can prove the theorem.
Let $v$ be a highest weight vector of $V$.
It is easy to see that in this case $Q_V^\dagger(\la)v=v$
(this follows from triangularity of $J$), and
$A_{w_0,V}(\la)v$ is constant.
This implies that $E_V(\la)v=v$.
But Part (i) of Lemma 2.15 in \cite{EV2} says that
a collection $E_V$ with the group-like property
such that $E_Vv=v$ on highest weight vectors
must necessarily be trivial: $E_V=1$.
The theorem is proved.

{\bf Remark.} Here is a proof of Corollary \ref{Aforsl(2)} using
Theorem \ref{Qdagger}, different from the proof in \cite{TV}. By
Theorem \ref{Qdagger}, $B^+(\la)=Q^\dagger(\la)$. Adapting the
formula of \cite{BBB} for the universal fusion matrix to our
conventions, we get that the element ${\mathcal J}(\la)$ has the
form
$$
{\mathcal J}(\la)=\sum_{k=0}^\infty
q^{-k(k+1)/2}\frac{(1-q^2)^k}{[k]_q!}
(f^n\otimes
e^n)\prod_{\nu=1}^n\frac{q^{2\la}}{1-q^{2\la+2\nu}(q^{-h}\otimes
  q^h)}.
$$
Applying the antipode to the first component, multiplying the
components, and changing $\la$ to $-\la-1$, we get the result.

\subsection{The determinant of $Q^\dagger$.}

Theorem \ref{Qdagger} allows us to
compute explicitly the determinant of $Q^\dagger$ on every
weight subspace of a finite dimensional module.

\begin{proposition} One has
$$
\text{det}(Q^\dagger(\la)|_{V[\beta]})=
\prod_{\alpha>0}
\prod_{k\ge 0}B^+_{\alpha\beta k}(q_\alpha,\la)^{d_V(\alpha,\beta,k)},
$$
where
$$
d_V(\alpha,\beta,k)=\dim(V[\beta+k\alpha])-
\dim(V[\beta+(k+1)\alpha]).
$$
\end{proposition}

\begin{proof} The proposition is immediate from Theorem
  \ref{Qdagger} and Proposition \ref{deter}.
\end{proof}

\section{Applications of the dynamical Weyl groups to trace
  functions}

In this section we will assume for simplicity that
$\g$ is a finite dimensional semisimple Lie
algebra, although the results
can be generalized to Kac-Moody algebras (see \cite{ES2}).

\subsection{A generalized Weyl character formula}

Recall that in \cite{EV2} we defined the trace functions
$$
\Psi^v(\la,\mu)= \text{Tr}|_{M_\mu}(\Phi_\mu^{v}
q^{2\la}),
$$
$$
\Psi_V(\la,\mu)=\sum \Psi^{v_i}(\la,\mu)
\otimes v_i^* \in V[0]\otimes V^*[0],
$$
where $v_i$ is a basis of $V[0]$ and $v_i^*$ is the dual basis of
$V^*[0]$.

Suppose $\mu$ is a large dominant integral weight, and $L_\mu$ is
the irreducible finite dimensional representation with this
highest weight. The intertwining operator $\Phi_\mu^v$
descends to an operator $\bar\Phi_\mu^v: L_\mu\to L_\mu \otimes V$,
and we define
$$
\Psi^v_\mu(\la)= \text{Tr}|_{L_\mu}(\bar\Phi_\mu^{v}
q^{2\la}),
$$
$$
\Psi_V^\mu(\la):=\sum \Psi^{v_i}_\mu(\la)\otimes v_i^*.
$$

Let us regard $\Psi_V(\la,\mu)$, $\Psi_V^\mu(\la)$ as linear
operators on $V[0]$.

\begin{proposition} \label{genWeyl}
One has
$$
\Psi_V^\mu(\la)=\sum_{w\in \W}(-1)^w \Psi_V(\la,w\cdot
\mu)A_{w,V}(\mu).
$$
\end{proposition}

\begin{proof} The proof is analogous to the proof of the
Weyl character formula using the approach of \cite{BGG};
it is based on the fact that in the Grothendieck group of
the category $\mathcal O$,
an irreducible module is an alternating sum of Verma modules.
\end{proof}

{\bf Remark 1.} If $V=\C$, this formula reduces to the Weyl
character formula.

{\bf Remark 2.} If $\g={{\frak sl}_n}$, $V=S^{rn}\C^n$
(the ``Macdonald'' case), then
$V[0]$ is 1-dimensional, and so the action of $A_{w,V}$ on $V[0]$
can be computed explicitly using the expression for $A_{s,V_m}$
for ${{\frak sl}_2}$. In this case it is easy to show that
Proposition \ref{genWeyl} reduces to Conjecture 8.2
in \cite{FV}, which was proved in \cite{ESt} (Prop. 5.3).

Recall that in \cite{EV2} we also defined renormalized trace
functions
$$
F_V(\la,\mu)=\delta_q(\la)\Psi_V(\la,-\mu-\rho)Q_V^{-1}(-\mu-\rho),
$$
where $\delta_q(\la)$ is the Weyl denominator:
$\delta_q(\la)=\sum_{w\in \W} (-1)^w q^{2(\la,w\rho)}$.

Let us say that a weight $\mu$ is antidominant if $-\mu$ is
dominant.
By analogy with the above, one can also define
for large antidominant integral weight $\mu$
$$
F_V^\mu(\la)=\delta_q(\la)\Psi_V^{-\mu-\rho}(\la)Q_V^{-1}(-\mu-\rho),
$$

\begin{corollary} \label{FWeyl} For an antidomonant $\mu$
with sufficiently large coordinates one has
$$
F_V^\mu(\la)=\sum_{w\in \W}(-1)^w F_V(\la,w\mu)(\A_{w,V^*}(\mu)^{-1})^*.
$$
\end{corollary}

\begin{proof} The proof follows from Proposition \ref{Adual},
and Proposition \ref{genWeyl}.
\end{proof}

\subsection{Dynamical Weyl group invariance of Macdonald-Ruijsenaars
operators and renormalized trace functions}

Let $q\ne 1$.
Recall that in \cite{EV2} we defined
the modified dynamical R-matrix $\RR_{UV}(\lambda)=R_{UV}(-\lambda-\rho)$,
and introduced Macdonald-Ruijsenaars operators, acting on
rational functions of $\lambda$ with values in $V[0]$
$$
{\mathcal D}_U=\sum_\nu
\text{Tr}|_{U[\nu]}(\RR_{UV}(\lambda))T_\nu,
$$
where $T_\nu$ maps $f(\lambda)$ to $f(\lambda+\nu)$.

\begin{proposition}\label{Macdonald} The Macdonald-Ruijsenaars operators
are invariant with respect to the dynamical action of $\W$ on
functions on $\h^*$ with values in $V[0]$. That is, $[{\mathcal
D_U}, w*]=0$ for $w\in \W$.
\end{proposition}

\begin{proof}
We have
$$
((w*)^{-1}{\mathcal D}_U(w*)f)(\la)=
\A_w(\la)^{-1}\sum_\nu \text{Tr}|_{U[w\nu]}(\RR(w\la))
\A_w(\la+\nu)f(\la+\nu)
$$
(for brevity we drop the subscripts indicating the modules in
which the operators act).
 From Corollary \ref{Rtrans} we easily obtain
$$
\RR(w\lambda)=
\A^{(2)}_{w}(\lambda)
\A^{(1)}_{w}(\lambda+h^{(2)})\RR_{VU}(\lambda)
\A^{(2)}_{w}(\lambda+h^{(1)})^{-1}
\A^{(1)}_{w}(\lambda)^{-1}
$$
Let us substitute this equation into the previous equation, and
use the fact that in the second component we are restricting to
the zero-weight subspace. It is easy to see that
the $\A$-factors cancel, and we get
$$
(w*)^{-1}{\mathcal D}_U(w*)={\mathcal D}_U,
$$
as desired.
\end{proof}

Let us return to the renormalized trace function $F_V(\la,\mu)$,
which we will now regard as an element of $V[0]\otimes V^*[0]$.

Recall from \cite{EV2} that $F_V(\lambda,\mu)$
 satisfies the Macdonald-Ruijsenaars equations
$$
{\mathcal
  D}_U^{(\lambda)}F_V(\lambda,\mu)=\chi_U(q^{-2\mu})F_V(\lambda,\mu),
$$
the dual Macdonald-Ruijsenaars equations
$$
{\mathcal
  D}_U^{(\mu)}F_V(\lambda,\mu)=\chi_U(q^{-2\la})F_V(\lambda,\mu),
$$
and has the symmetry property $F_V(\la,\mu)=F_{V^*}(\mu,\la)^*$
(here the superscripts $\la,\mu$ denote the variables
with respect to which to take shifts when applying difference operators,
and $()^*$ denotes the operator of exchanging the factors $V[0]$ and
$V^*[0]$).

\begin{proposition}\label{tracefun} The function $F_V(\lambda,\mu)$
is invariant under the dynamical action of $\W$
on functions with values in $V[0]$. That is,
$$
F_V(\lambda,\mu)=(\A_{w,V}(w^{-1}\lambda)\otimes
\A_{w,V^*}(w^{-1}\mu)) F_V(w^{-1}\lambda,w^{-1}\mu),\ w\in \W.
$$
\end{proposition}

\begin{proof} The proof is similar to the proof
of the symmetry of $F_V$, given in \cite{EV2}.

It suffices to assume that $q\ne 1$. Let $F_V'(\la,\mu)$ denote
the right hand side of the equality to be proved. By Proposition
\ref{Macdonald}, $F_V'$, like $F_V$, is a solution of the
Macdonald-Ruijsenaars equations and the dual
Macdonald-Ruijsenaars equations. Moreover, both $F_V$, $F_V'$
have the form: $q^{2(\la,\mu)}$ times a finite sum of rational
functions of $q^{(\la,\alpha_i)}$ multiplied by rational functions of 
$q^{(\mu,\alpha_i)}$
(where the denominators of the rational functions are products of
binomials of the form $1-q^{(\la,\beta)}$, respectively $1-q^{(\mu,\beta)}$).

Let us regard $F_V,F_V'$ as functions with values in
$\text{End}(V[0])$. It is easy to see, using power series
expansions, that a solution of the Macdonald-Ruijsenaars equations
with the above properties is unique up to right multiplication by
an operator depending rationally of $q^{(\mu,\alpha_i)}$. Similarly, a solution
of the dual Macdonald-Ruijsenaars equations with such properties
is unique up to left multiplication by an operator depending
rationally on $q^{(\la,\alpha_i)}$. So we have
$$
F_V'(\la,\mu)=X(\la)F_V(\la,\mu),\qquad
F_V'(\la,\mu)=F_V(\la,\mu)Y(\mu),
$$
where $X,Y$ are rational operator valued functions of
$q^{(\la,\alpha_i)}, q^{(\mu,\alpha_i)}$, and hence
$$
X(\la)F_V(\la,\mu)=F_V(\la,\mu)Y(\mu).
$$

Let us take the limit $q^{(\la,\alpha_i)}\to 0$ (for all $i$) in
the last equality. It follows from the asymptotics of
intertwiners (see \cite{ESt}) that in this limit $F_V$ is
equivalent to $q^{-(\la,\mu)}\text{Id}$. So we get $\lim
X(\la)=Y(\mu)$ for all $\mu$. Thus, $Y(\mu)$ is a constant
operator. Using the symmetry of $F_V$, we get that $X(\la)$ is
also a constant, so we get $X(\la)=Y(\mu)=X$, where $X$ is a
constant operator.

Finally, let us show that $X=1$. We have the identity
$$
XF_V(\la,\mu)=\A_{w,V}(w^{-1}\la)F_V(w^{-1}\la,w^{-1}\mu)
\A_{w,V^*}(w^{-1}\mu)^*.
$$
Using Proposition \ref{dualonweight0}, we can rewrite this
equation in the form
$$
XF_V(\la,\mu)=\A_{w,V}(w^{-1}\la)F_V(w^{-1}\la,w^{-1}\mu)
\A_{w^{-1},V}(-\mu).
$$
Now let us take the limit: $q^{(\mu,\alpha_i)}\to 0$,
$q^{(w^{-1}\la, \alpha_i)}\to 0$. Then, using Corollary
\ref{limits1.2}, we get
$$
X=A_{w,V}^-A_{w^{-1},V}^+=1,
$$
as desired.
\end{proof}

\subsection{The multicomponent dynamical action
and invariance of multicomponent trace functions}

Let $V_1,...,V_N$ be integrable $U_q(\g)$-modules.
Define the linear operator
$A_{w,V_1,...,V_n}(\la): V_1\otimes...\otimes V_N\to
V_1\otimes...\otimes V_N$ by
$$
A_{w,V_1,...,V_N}(\la)=A^{(N)}_{w,V_N}(\la)
A^{(N-1)}_{w,V_{N-1}}(\la-h^{(N)})...
A^{(1)}_{w,V_1}(\la-h^{(2)}-...-h^{(N)}).
$$
With these operators one can associate the action
of $\tilde \W$ on functions of $\la$ with
values in $V_1\otimes...\otimes V_N$ given by
$$
(w\bullet f)(\la)=A_{w,V_1,...,V_N}(w^{-1}\cdot \la)
f(w^{-1}\cdot \la).
$$
We call this action the {\bf shifted multicomponent
dynamical action}. As before, the (unshifted)
{\bf multicomponent dynamical action} is defined by
$$
w\diamond f=\A_{w,V_1,...,V_N}(w^{-1} \la)
f(w^{-1}\la),
$$
where
$$
\A_{w,V_1,...,V_N}(\la)=
A_{w,V_1,...,V_N}\left(-\la-\rho+\frac{1}{2}(\sum_{i=1}^N h^{(i)})\right)
$$

Recall from \cite{EV2} that the operator $J^{1...N}(\la)$
on $V_1\otimes...\otimes V_N$ is defined by
$$
J^{1...N}(\la)=J^{1,2...N}(\la)...J^{N-1,N}(\la).
$$

\begin{lemma} \label{multi}
Conjugation with
the operator $J^{1...N}$ transforms the shifted dynamical action
of $\tilde \W$ into its shifted multicomponent dynamical action.
That is,
$$
J^{1...N} (w\bullet)= (w\circ)J^{1...N}.
$$
\end{lemma}

\begin{proof}
This follows by applying Lemma \ref{com-A}
several times.
\end{proof}

Recall from \cite{EV2} the definition of the quantum KZB
operators $K_j,K_j^\vee$ and the diagonal operators $D_j,D_j^\vee$
acting on functions of
 $\la$ and $\mu$ with values in $(V_1\otimes...\otimes V_N)[0]$, respectively
 $(V_N^*\otimes...\otimes V_1^*)[0]$ ($j=1,...,N$). Namely, set
 $$
 D_j=(q^{-2\mu-\sum x_i^2})_{*j}(q^{-2\sum x_i\otimes x_i})_{*j,*1...*j-1},
 $$
 (where $*i$ labels the component $V_i^*$), and
 $$
 K_j=
 \RR _{j+1,j}(\la+h^{(j+2...N)})^{-1}...\RR _{Nj}(\la)^{-1}
 \Gamma_j\times
$$
$$
\RR _{j1}(\la+h^{(2...j-1)}+h^{(j+1...N)})...\RR _{jj-1}
 (\la+h^{(j+1...N)}),
 $$
 where $\Gamma_jf(\la):=f(\la+h^{(j)})$, and $h^{j...k}$
 acting on a
 homogeneous multivector has to be replaced with the sum of weights of
 components $j,...,k$ of this multivector.
 Analogously, define the operators
 $$
 D_j^\vee=
 (q^{-2\la-\sum x_i^2})_{j}(q^{-2\sum x_i\otimes x_i})_{j,j+1...N},
 $$
 and
 $$
 K_j^\vee=
 \RR _{*j-1,*j}(\mu+h^{(*1...*j-2)})^{-1}...\RR _{*1,*j}(\mu)^{-1}
 \Gamma_{*j}\times
$$
$$
 \RR _{*j,*N}(\mu+h^{(*j+1...*N-1)}+h^{(*1...*j-1)})...
 \RR _{*j,*j+1}
 (\mu+h^{(*1...*j-1)}),
 $$
 where $\Gamma_{*j}f(\mu)=f(\mu+h^{*j})$.

In \cite{EV2}, we defined
the multicomponent
renormalized trace functions $F_{V_1,...,V_N}(\la,\mu)$
with values in $(V_1\otimes...\otimes V_N)[0]\otimes
(V_N^*\otimes...\otimes V_1^*)[0]$. Two of
our main results were the identities
$$
(K_j\otimes D_j)F_{V_1,...,V_N}=F_{V_1,...,V_N}
$$
(the quantum KZB equations), and
$$
(D_j^\vee\otimes K_j^\vee)F_{V_1,...,V_N}=F_{V_1,...,V_N}
$$
(the dual quantum KZB equations).

\begin{corollary} (i) The multicomponent
renormalized trace functions $F_{V_1,...,V_N}(\la,\mu)$
are invariant under the multicomponent
dynamical action $w\diamond$ of $\W$ in both components. That is:
$$
((w\diamond)^{(\la)}\otimes (w\diamond)^{(\mu)})
F_{V_1,...,V_N}(\la,\mu)=F_{V_1,...,V_N}(\la,\mu).
$$

(ii) The quantum KZB operators $K_j,K_j^\vee$ and
the diagonal operators $D_j,D_j^\vee$
are invariant under the multicomponent dynamical action of
$\W$. In particular, the qKZB
and dual qKZB equations are invariant
under the multicomponent dynamical action.
\end{corollary}

\begin{proof} Statement (i) follows from Lemma \ref{multi},
Proposition \ref{tracefun} and the definitions of \cite{EV2}.
Statement (ii) can be
checked directly using Corollary \ref{Rtrans}.
\end{proof}

\section{Dynamical Weyl groups for affine Lie algebras
and quantum affine algebras}

In this section we will consider the
dynamical Weyl group in the case when
the role of $\g$ is played by an affine Kac-Moody Lie algebra
$\tilde \g$,
and the role of integrable representations $V$ of $\g$ or $U_q(\g)$ is
played by representations on Laurent polynomials in one variable
with coefficients in a finite dimensional vector space
(we call them loop representations).

This situation
turns out to be especially interesting for applications.
Although this setting is very similar
to the one already considered, it is not exactly the same, since
loop representations are not integrable.
Therefore, we will describe the changes that are
necessary to carry out our main construction in this new situation.

\subsection{Affine algebras
and loop representations}

In this subsection we will recall 
some standard facts about finite dimensional representations of classical and 
quantum affine algebras. The material on the classical affine algebras 
is standard; most of the material on the quantum case can be found in 
the book \cite{CP}, and references therein.

Let $\g$ be a simple finite dimensional Lie algebra.
Let $(,)$ be the form on $\g$ defined in Section 2, and
let the positive integer  $m$ be defined by
$(\theta,\theta)=2m$ for the maximal root $\theta$ of $\g$ (we
have $m=1$ in the simply laced case, $m=3$ for $\g=G_2$,
and $m=2$ otherwise).
Let $\hat \g=\g[z,z^{-1}]\oplus \C c$
be the standard central extension of the loop algebra:
$$
[a(z),b(z)]=[ab](z)+m\text{Res}_0(a'(z),b(z))c.
$$
Let $\tilde \g=\C d\oplus \hat \g$ be the extension
of $\hat\g$ by the scaling element $d$ such that
$[d,a(z)]=za'(z)$, $[d,c]=0$.

The Lie algebra $\tilde\g$ is the Kac-Moody Lie algebra
corresponding to the affine (i.e. extended) Cartan matrix of $\g$.
In particular, we have $\tilde\g=\tilde\h\oplus \hat \n_+\oplus \hat \n_-$,
where $\hat \n_\pm$ are the positive and
negative nilpotent subalgebras,
$\tilde \h=\h\oplus \C c\oplus \C d$.
The restriction of the invariant bilinear form on $\tilde\g$
to $\tilde\h$ is defined by $(d,d)=(c,c)=0$, $(d,c)=1/m$,
$(d,h)=(c,h)=0$, $h\in \h$ (and the form on $\h$ is the same as
in $\g$).

{\bf Remark.} This normalization of the invariant form is
traditional in the theory of quantum groups. On the other hand,
the traditional bilinear form in the representation theory of
affine Lie algebras and KZ equations is defined by the condition
$(\theta,\theta)=2$ on the dual space, so it is $m$ times bigger
than our form. This is why many of our formulas have an extra
factor $m$ compared to the formulas from other texts about KZ
equations (e.g. \cite{EFK}).

The dual Cartan subalgebra of $\tilde\g$ can be written as $\tilde
\h^*=\h^*\oplus \C c^*\oplus \C d^*$, where $c^*,d^*$ are the
dual elements to $c,d$. Thus, elements of $\tilde\h^*$ can be
written as triples $(\la,k,\Delta)$, where $\la\in \h^*$,
$k,\Delta\in \C$, i.e. $(\la,k,\Delta)(h+ac+bd)=\la(h)+ka+\Delta
b$. The number $k$ is called the central charge of $\tilde\la$.
For instance, the roots of $\tilde\g$ are the elements of the
form $(\alpha,0,n)$, where $n\in \Z$, $\alpha$ is $0$ or is a
root of $\g$, and $(\alpha,n)\ne (0,0)$. The special weight
$\rho_{\tilde \g}$ for the affine algebra $\tilde\g$ will be
denoted by $\tilde\rho$. It has the form $\tilde\rho=\rho+h^\vee
c^*$, where $h^\vee$ is the dual Coxeter number of $\g$, and
$\rho$ is the special weight for the finite dimensional Lie
algebra $\g$, regarded as an element of $\tilde \h^*$. In other
words, $\tilde\rho=(\rho,h^\vee,0)$.

We will denote the Chevalley generators of $\g$ by $e_i,f_i,h_i$,
$i=1,...,r$, and the additional generators of $\hat\g$ by $e_0,f_0,h_0$.
Similarly, $\alpha_1,...,\alpha_r$ will stand for simple roots of
$\g$, and $\alpha_0$ for the additional simple root of $\hat\g$.

Let $U_q(\tilde\g)$, $U_q(\hat \g)$,
$U_q(\hat \n_\pm)$ be the quantum deformations
of the corresponding classical objects, defined as in the
general Kac-Moody case. In particular,
the algebra $U_q(\tilde\g)$ contains elements $q^{bc}$, $q^{bd}$,
$b\in \C$.

As before, we will consider the quantum situation but will allow
$q$ to be $1$, unless otherwise specified.

Let us define the notion of a
{\bf loop representation}
of $U_q(\tilde\g)$. Let $\bar V$ be a
finite dimensional representation
of $U_q(\hat\g)$. Set $V=\bar V[z,z^{-1}]$, with the following
action of $U_q(\tilde\g)$:
$d|_{V}=z\frac{d}{dz}$, and
$x|_{V}=z^nx|_{\bar V}$ for $x\in U_q(\hat\g)$,
such that $[d,x]=nx$.

\begin{definition} $V$ is called a loop
  representation.
\end{definition}

Thus every loop representation $V$ has a
natural structure of a module over $\C[z,z^{-1}]$, and
the underlying representation $\bar V$ is reconstructed by $\bar
V=V/(z-1)V$. Moreover, if $a\in \C^*$ then we get
a new finite dimensional representation $\bar V(a)=V/(z-a)V$ of
$U_q(\hat\g)$.
This representation is called the shift of $\bar V$ by $a$.

We will need the following proposition.

\begin{proposition}\label{cis0} (i) Any finite dimensional representation
$Y$ of $U_q(\hat \g)$ has a weight decomposition
with respect to $U_q(\hat\h)\subset U_q(\hat\g)$
(where $\hat\h:=\h\oplus \C c$).

(ii) In any finite dimensional representation $Y$
of $U_q(\hat\g)$, the element $c$
acts by  zero (in the q-case, by this we mean that $q^{bc}=1$).

(iii) Statements (i) and (ii) are valid for loop
representations.
\end{proposition}

\begin{proof} This Proposition is well known, but we will give
a proof for the reader's convenience.

 Statement (i) follows from
the fact that any finite dimensional representation has a weight
decomposition with respect to any (quantum) ${{\frak
sl}_2}$-subalgebra corresponding to a simple root (by
representation theory of quantum ${{\frak sl}_2}$).

Let us prove (ii). By the existence of a weight decomposition, it
suffices to prove this for irreducible representations. But in an
irreducible representation, $c$ (respectively, $q^{bc}$) acts by
a scalar. If $q=1$, this scalar must be zero, as $c$ is a linear
combination of $[e_i,f_i]$ and thus $\text{Tr}(c)=0$. If $q\ne
1$, it suffices to show that $q^c=1$ (as the weights are
integral). But
$(q_0-q_0^{-1})[e_0,f_0]=q^{c-\theta^\vee}-q^{-c+\theta^\vee}$.
Since $w_0\theta^\vee=-\theta^\vee$, in a finite dimensional
$U_q(\g)$-module we have $Tr(q^{\theta^\vee})=
Tr(q^{-\theta^\vee})\ne 0$. Thus,
$$
0=(q_0-q_0^{-1})Tr([e_0,f_0])=(q^c-q^{-c})Tr(q^{\theta^\vee}).
$$
Thus, $q^c=\pm 1$, but it is an integer power of $q$, so
$q^c=1$ as desired.

Statement (iii) is clear from (i),(ii).
\end{proof}

The main examples of finite dimensional representations of
$U_q(\hat\g)$ are the so called evaluation modules.
To define them, let us first assume that $q=1$.
In this case, for any $a\in \C^*$,
we have the evaluation homomorphism $\text{ev}_a: U(\hat\g)\to
U(\g)$, defined by $\text{ev}_a(x(z))=x(a),\text{ev}_a(c)=0$.
Let $Y$ be a finite dimensional irreducible $\g$-module.
Then let $Y(a)$ denote the $\hat\g$-module
$\text{ev}_a^*Y$ (the pullback of $Y$).
This module is called {\it an evaluation module}.
It is easy to see that the associated loop representation is
$Y[z,z^{-1}]$, with pointwise action of $\hat\g$.

For $q\ne 1$, by evaluation modules over $U_q(\hat \g)$
we will mean q-deformations of evaluation modules for
$\hat\g$; in other words, finite dimensional modules which remain
irreducible when restricting to $U_q(\g)$.

{\bf Remark.} For $\g={\frak sl}_n$, there exists
 an analog of the homomorphism $\text{ev}_a$,
introduced by Jimbo (see e.g. \cite{EFK}). In this case, we can define
the evaluation module $Y(a)$ corresponding to any irreducible
finite dimensional $U_q(\g)$-module $Y$, in the same way as in
the classical case. In other words, we can q-deform every
evaluation module. Outside of type A, the evaluation homomorphism
does not exist, and, as a result, not every evaluation module can
be deformed (e.g. the module corresponding to the adjoint
representation of $\g$ cannot); but some evaluation modules can still be
deformed, e.g. the vector representation for classical groups.

\subsection{The affine Weyl group}

In this subsection we will recall basic facts about affine Weyl
groups. These facts are standard, and can be found
in the literature (e.g. \cite{Ch1,Ch2} and references therein),
but we will give
the definitions, statements, and even some proofs
for the reader's convenience.

Let $\W^a$ denote the Weyl group of $\tilde\g$.
It has generators $s_0,...,s_r$, satisfying the usual braid and
involutivity relations.

Let $Q^\vee$ be the dual root lattice of $\g$. It is well known
that $\W^a$ is isomorphic to the semidirect product $\W \ltimes
Q^\vee$ of the Weyl group $\W$ of $\g$ with the dual root lattice
$Q^\vee$ (i.e. the Cartesian product $\W\times Q^\vee$ with the
product $(w,q)(w',q')=(ww',(w')^{-1}(q)+q')$), via the isomorphism
defined by $s_i\to (s_i,0), i\ne 0; s_0\to (s_\theta,
-\theta^\vee)$. In particular, $\W$ and $Q^\vee$ are subgroups of
$\W^a$ in a natural way. To avoid confusion, given an element
$\beta\in Q^\vee$, we will write $t_\beta$ for the corresponding
element of $\W^a$, and use the multiplicative (rather than the
additive) notation for the product of such elements; thus, $t_\mu
t_\nu=t_{\mu+\nu}$.

Let us compute the action of $\W^a$
on $\tilde\h^*=\h^*\oplus \C c^*\oplus \C d^*$.
The action of $\W$ is obvious
(i.e. it acts only on the $\h^*$-component, keeping the other two
unchanged), so let us calculate
the action of the lattice $Q^\vee$.

Recall that elements of $\tilde\h^*$ can be written as triples
$(\la,k,\Delta)$.

We have the following result.

\begin{lemma}\label{actonweights}
One has
$$
t_\nu(\la,k,\Delta)=(\la+mk\nu,k,\Delta-(\la,\nu)-mk(\nu,\nu)/2).
$$
\end{lemma}

\begin{proof} Let us call the operator defined by the right hand
  side by $t_\nu'$.
Using the identity $s_0s_\theta=t_{\theta^\vee}$,
we get that $t_\nu=t_\nu'$ for $\nu=\theta^\vee$.
Since the statement that $t_\nu=t_\nu'$ is Weyl group invariant,
and $t_{\nu_1+\nu_2}'=t_{\nu_1}'t_{\nu_2}'$, this is sufficient.
\end{proof}

It is also useful to introduce the extended affine Weyl group
$\W^b$. By the definition, $\W^b=\W\ltimes P^\vee$, where
$P^\vee$ is the dual weight lattice. Thus,
$\W^b$ naturally contains $\W^a$ as a subgroup.
We can define the action of $\W^b$ on weights
by extending the formula of Lemma \ref{actonweights}
to elements $\nu\in P^\vee$. It is easy to check that the set
of roots is invariant under $\W^b$ (as the pairing
between $P^\vee$ and the root lattice $Q$ takes only integer
values: these two lattices are dual to each other).

Let $G$ be the simply connected Lie group corresponding to $\g$.
It is easy to see that the exponential map
$\varepsilon:=\exp (2\pi i\cdot *): \g\to G$
identifies the group $P^\vee/Q^\vee=\Pi$ with
the center $G$. Indeed, the lattice $Q^\vee$ is the
kernel of $\varepsilon$ restricted to $\h$, so we have an
injective induced map
$\varepsilon: P^\vee/Q^\vee\to G$.
This map lands in the center since elements of $P^\vee$
give integer inner products with roots, and hence
elements $\varepsilon(x)$, $x\in P^\vee/Q^\vee$ act
by the same scalar on all weight subspaces of any irreducible
finite dimensional $G$-module. Reversing this argument, we see
that this map is also surjective, so it is an isomorphism.

This implies that
$\W$ acts trivially on $\Pi$. Indeed,
the action of $\W$ on $\Pi$ is induced by the action on the
maximal torus $T\subset G$ of the normalizer $NT\subset G$ of
this torus by conjugation. So the elements of $\Pi$ are
invariant under this action because they are central in $G$.

Thus, the subgroup $\W^a$ is normal in $\W^b$, and the quotient
$\W^b/\W^a$ is naturally identified with $\Pi$ (the
identification is induced by the embedding $P^\vee\to \W^b$).

It is useful to introduce the notion of the length
of an element of $\W^b$. By the definition, let the length
of $w\in \W^b$, denoted by $l(w)$, be the number of positive roots
which are mapped under $w$ to negative roots.
It is obvious that this number is finite,
and that simple reflections have length 1.
It is known
 that the length of an element of $\W^a$ given by
a reduced decomposition with $n$ factors is $n$, and that if
$\la,\mu\in P^\vee$ are dominant then
$l(\la+\mu)=l(\la)+l(\mu)$ (see e.g. \cite{Ch1,Ch2}, and references therein; 
the statements
follow by looking at how $\la$ and $\mu$ act on positive roots).

Let $\tilde\Pi\subset \W^b$ be the group of transformations
that have length $0$, i.e. those which map the sets of positive and negative
affine roots to themselves.

It is clear that any element $w\in \W^b$ of length $n>0$ can be
represented as a product $w=\sigma w'$, where $w'$ has length
$n-1$. Indeed, let $\alpha$ be a simple positive root
such that $w^{-1}\alpha$ is negative (clearly such exists,
otherwise
$w\in\tilde \Pi$).
If $\beta$ and $w\beta$ are positive roots,
then $w\beta\ne \alpha$, so $s_\alpha
w\beta$ is positive.
In addition, $s_\alpha w (-w^{-1}\alpha)=\alpha$ is positive, so
$l(s_\alpha w)=n-1$.

Thus, we get a factorization $\W^b=\tilde\Pi\W^a$. Moreover, the
factorization  of an element of $\W_b$ into a product of elements
of $\W_a$ and $\tilde \Pi$ is unique, since $\W^a$ and $\tilde
\Pi$ intersect trivially (as nontrivial elements of $\W^a$ have
positive length). In other words, the exact sequence
$$
1\to \W^a\to \W^b\to \Pi \to 1
$$
is split (canonically!), and the subgroup $\W^a\subset \W^b$ is complemented by
the subgroup $\tilde\Pi$. Thus, we have $\W^b=\tilde\Pi
\ltimes \W^a$. 

In fact, the canonical splitting homomorphism
$\eta: \Pi=P^\vee/Q^\vee\to \W^b$ can be constructed
explicitly as follows.

For any $i=1,...,r$, let $w_0^i$ be the maximal element of the
Weyl group of the Levi subalgebra of $\g$,
whose Dynkin diagram is obtained from that of $\g$
by throwing away the i-th vertex of the Dynkin diagram.
Let $w_{[i]}=w_0w_0^i$.

Recall that the fundamental coweights for $\g$,
$\omega_i^\vee$, $i=1,...,r$, are the elements of $P^\vee$
defined by $\alpha_j(\omega_i^\vee)=\delta_{ij}$.
We say that a fundamental coweight $\omega_i^\vee$ is minuscule if
$\theta(\omega_i^\vee)=1$. (For non-minuscule
coweights, this number is greater
than 1).

It is known that for any element $\pi\ne 1$ of $\Pi$,
there exists a unique minuscule fundamental coweight
$\omega_i^\vee\in P^\vee$ which represents $\pi$ in $P^\vee/Q^\vee$, and
all minuscule fundamental coweights are
obtained in this way exactly once.

Let us denote the coweight $\omega_i^\vee$ corresponding to $\pi$
by $\omega^\vee(\pi)$, and the index $i$ by $i_\pi$.

\begin{proposition}
The homomorphism $\eta$ is defined by
$\eta(\pi)=t_{\omega^\vee(\pi)}w_{[i_\pi]}^{-1}$.
\end{proposition}

\begin{proof} We need to show that
$\eta$ lands in $\tilde\Pi$ and that it is a homomorphism.
Let us prove the first statement. So let $\pi\in \Pi$,
$i_\pi=i$, and let us show that $\eta(\pi)\in \tilde\Pi$.

It is clear
that if $j\ne 0,i^*$ (where $i^*$ is the dual vertex to $i$)
then $\alpha_j$ is mapped under $\eta(\pi)$ to
$\alpha(j')$, with $j'\ne 0$. Thus, we need to show that
$\alpha_{i^*}$ and $\alpha_0$ are also mapped to simple positive
roots. A simple computation shows that this property is
equivalent to the identity $w_0^i\alpha_i=\theta$. So let us
prove this.

Let $w_0^i\alpha_i=\beta$. Clearly, $(\beta,\omega_i^\vee)=1$.
So $\beta$ is a positive root of the form
$\beta=\alpha_i+\gamma$, where
$\gamma$ is a linear combination of positive
roots except $\alpha_i$. Similarly,
$\theta=\alpha_i+\gamma'$, and $\gamma'\ge \gamma$.
Thus,
$$
w_0^i\theta=w_0^i\alpha_i+w_0^i\gamma'=
\alpha_i+\gamma+w_0^i\gamma'
$$
Now, we see that since $\gamma'\ge \gamma$,
the height of the right hand side (i.e. the
sum of the multiplicities of the simple roots) is $\le 1$,
and the equality is possible only if $\gamma=\gamma'$.
But the right hand side is a positive root,
so $\gamma'=\gamma$ and hence $\beta=\theta$.

Now we prove the second statement (that $\eta$ is a homomorphism).
Since $\eta$ lands in $\tilde\Pi$, it is sufficient to check
that the map $\Pi\to \W^b/\W^a=P^\vee/Q^\vee$
induced by $\eta$ is a homomorphism. But this is obvious from the
definition.
\end{proof}

Thus, $\Pi$ can be identified with $\tilde\Pi$
via $\eta$, and can thus be regarded as a subgroup of $\W^b$;
we will assume from now on that we have performed this identification.
 In particular, any element $\pi\in \Pi$ acts on
$\W^a$ by conjugation. It is easy to show that this action
permutes simple reflections according to an automorphism of the
extended Dynkin diagram $\Gamma^a$ of $\g$. In other words, we
have a homomorphism (in fact, an embedding) $\Pi\to
\text{Aut}(\Gamma^a)$.

{\bf Remark.} It is easy to check that $\text{Aut}(\Gamma^a)=
\text{Aut}(\Gamma)\ltimes \Pi$, where $\Gamma$ is the Dynkin
diagram of $\g$.

{\bf Example.} Consider $\g={\frak sl}_2$. In this case the group
$\W^a$ is generated by two elements $s_0,s_1$ such that
$s_0^2=s_1^2=1$, with no other relations. So we can think of
$\W^a$ as the group of all affine linear transformations 
of the integers, which is
generated by $s_0(x)=-x$, $s_1(x)=1-x$. The group $\W^b$ is in
this case the set of all affine linear transformations 
of the half-integers, so it has the form
$\W^b=\Pi\ltimes \W^a$, where $\Pi=\lbrace{1,\pi\rbrace}$,
$\pi(x)=\frac{1}{2}-x$. We see that
$Q^\vee=\Z,P^\vee=\frac{1}{2}\Z$, and the action of $\Pi$ on
$\W^a$ is given by $\pi s_0\pi^{-1}=s_1, \pi s_1\pi^{-1}=s_0$, as
predicted by the general theory. We have $\omega_1^\vee(x)=x+1/2$,
and the element $w_{[1]}$ is given by $w_{[1]}(x)=-x$.

\subsection{Intertwining operators and expectation values}

Let $X$ be a module over $\tilde\h$ which has a weight
decomposition, and let $M$ be a module over $U_q(\tilde\g)$ from
category
$\mathcal O$.
For any weight $\tilde\nu$, define the space
$(M\hat\otimes X)[\tilde\nu]$ to be
$\hat\oplus_{\tilde\beta}M[\tilde\beta]\otimes X[\tilde\nu-\tilde\beta]$
(where $\hat\oplus$ is the completed direct sum, i.e. the Cartesian
product over all $\tilde\beta$). Elements of this space are arbitrary
(possibly infinite) sums of tensors whose first and second
components are homogeneous.
Define the completed tensor product $M\hat\otimes X$ to be
$\oplus_{\tilde\nu} (M\hat\otimes X)[\tilde\nu]$ (an algebraic direct sum).

Let $M_{\tilde\la}$ be the Verma module over $U_q(\tilde\g)$
with highest weight $\tilde\la$.
Let $\bar V$ be a finite dimensional representation of $U_q(\hat
\g)$, and $V$ the corresponding loop representation.
Consider an intertwining operator
$\Phi: M_{\tilde\lambda}\to M_{\tilde\mu}\hat\otimes V$.
Like for intertwiners into the usual tensor product,
we have $\Phi v_{\tilde\lambda}=v_{\tilde\mu}\otimes v+...$, where $...$
denote terms of lower weight in the first component, and
$v\in V[\tilde\lambda-\tilde\mu]$ (but now the sum
denoted by ... may be infinite).
By analogy with the previous setting,
we will call $v$ the expectation value of
$\Phi$ and write $<\Phi>=v$.

Let us now generalize Lemma \ref{inter} to the affine case.
For clarity we will split this generalization into two lemmas.

Let $V$ be a loop representation, and $\tilde\nu$ a weight in
$V$.

\begin{lemma}\label{inter1} For generic $\tilde\lambda$
the map $\Phi\to <\Phi>$ is an isomorphism of
  vector spaces
$\text{Hom}_{U_q(\tilde\g)}(M_{\tilde\lambda},M_{\tilde\lambda-\tilde\nu}
\hat\otimes V)\to
V[\tilde\nu]$.
\end{lemma}

{\bf Remark 1.} Here ``generic'' means away from a countable
(possibly infinite) number
of hyperplanes.

{\bf Remark 2.} Note that the lemma would be wrong
if we used $\otimes$ instead of $\hat\otimes$.

{\bf Remark 3.} Note that in Lemma \ref{inter1}, the central
charges of $\tilde\la$ and $\tilde\la-\tilde\nu$ are the same, by
Proposition \ref{cis0}, unless the spaces are zero.

\begin{proof} Straightforward,
as in Lemma \ref{inter}; see also Theorem 3.1.1 in \cite{EFK}.
\end{proof}

This lemma allows one to define the interwining
operator $\Phi_{\tilde\la}^v$ with expectation value $v$.
As before, it has coefficients which are rational functions of
$\tilde\la$
or $q^{(\tilde\la,\alpha_i)}$.

\begin{lemma}\label{inter2}
The map $\Phi\to <\Phi>$ is an isomorphism for dominant weights
$\tilde\la$ with sufficiently large coordinates $\tilde\la(h_i)$.
\end{lemma}

\begin{proof} The lemma is proved by arguments
similar to those in \cite{ESt}. Namely,
similarly to \cite{ESt}, one can write down an explicit formula
for $\Phi_{\tilde\la}^v v_{\tilde\la}$,
and show that its poles are all of first
order and can occur only on hyperplanes $(\tilde\la+\tilde\rho,\alpha)=
\frac{n}{2}(\alpha,\alpha)$ for positive roots $\alpha$ and such
$n>0$ that $V[\tilde\nu+n\alpha]\ne 0$. If a dominant weight $\tilde\la$
belongs to such a hyperplane, then $(\alpha,\alpha)>0$, so
$\alpha$ is a real root. But it is clear that there exists a number
$N$ such that for $n\ge N$ one has $V[\tilde\nu+n\alpha]=0$ for any
weight $\tilde\nu$ of $V$
and any real root $\alpha$.
\end{proof}

{\bf Remark.} Note that the last statement of the proof of Lemma
\ref{inter2}
would be false for imaginary roots.

\subsection{The dynamical Weyl group for
loop representations}

It is easy to see that loop representations
are locally finite, so the dynamical Weyl group operators
$A_{w,V}(\tilde\la)$, $w\in \W^a$, are already defined on them.
It is easy to see that these operators are linear
over $\C[z,z^{-1}]$.

Moreover, we have
an analog of Proposition \ref{restr}:

\begin{proposition}\label{restr1}  One has
$$
\Phi v_{w\cdot \tilde\la}^{\tilde\la}=v_{w\cdot \tilde\mu}^
{\tilde\mu}\otimes
A_{w,V}(\tilde\la)<\Phi>+\text{lower weight
  terms}.
$$
\end{proposition}

\begin{proof} The proof is analogous to the proof of
Propsition \ref{restr}.
\end{proof}

\subsection{Fusion matrices}

Now let us generalize to the affine case the construction of
fusion matrices.

First of all, define completed tensor products
of Laurent polynomial spaces. Let
$\bar V_i$ be vector spaces,
and $V_i=\bar V_i[z,z^{-1}]$.
Define $V_1\ovo ...\ovo V_N$ to be
$$
V_1\ovo ...\ovo V_N:=
(\bar V_1\otimes...\otimes
\bar V_n)[[z_2/z_1,...,z_N/z_{N-1}]][z_1,z_1^{-1}],
$$
where $z_i$ denote the formal parameters corresponding to $V_i$.
It is clear that if $V_i$ are loop representations
of $U_q(\tilde\g)$ then $V_1\ovo ...\ovo V_N$ is also a
representation of this algebra, which is locally finite.
In fact, $V_1\ovo...\ovo V_N$ is a certain completion
of the ordinary tensor product $V_1\otimes...\otimes V_N$.

Now
let $\tilde\la\in \tilde\h^*$ be a generic weight.
Let $V,U$ be loop representations of $U_q(\tilde\g)$, and
$v \in V[\tilde\mu],\;u\in U[\tilde\nu]$.

Consider the composition
\bean
\Phi^{u,v}_{\tilde\lambda}:\;M_{\tilde\lambda}
\stackrel{\Phi^v_{\tilde\lambda}}{\longrightarrow}
M_{\tilde\lambda-
\tilde\mu}
\hat\otimes V \stackrel{\Phi^u_{\tilde\lambda-\tilde\mu}\otimes 1}
{\longrightarrow}
M_{\tilde\lambda-\tilde\mu-\tilde\nu} \hat\otimes (U \ovo V).
\notag
\eean
(It is easy to see that this compoosition is well defined; 
see also \cite{EFK} for explanations). 
Then
$\Phi^{u,v}_{\tilde\lambda} \in
\mathrm{Hom}_{U_q(\g)}(M_{\tilde\lambda},M_{\tilde\lambda-\tilde\mu-\tilde\nu}
\hat\otimes (U \ovo V))$. Let $x=<\Phi^{u,v}_{\tilde\lambda}>$.
Since $U,V$ have a weight decomposition by Proposition \ref{cis0},
the assignment
$(u, v) \mapsto x$ is bilinear, and naturally extends to
a zero weight map
$$
J_{UV}(\tilde\lambda):\; U \ovo V \to U \ovo V,
$$
linear over $\C[[z_2/z_1]][z_1,z_1^{-1}]$.
This means, the operator $J_{UV}(\tilde\la)$ can be understood as
an element of $\text{End}(\bar U\otimes \bar V)[[z_2/z_1]]$.

The operator $J_{UV}(\tilde\lambda)$
is called the fusion matrix of $U$ and $V$.
The fusion matrix $J_{UV}(\tilde\lambda)$ is a
power series in $z=z_2/z_1$ of the form
$$
J_{UV}(\tilde\la)(z)=\sum_{n\ge 0} J_{UV,n}(\tilde\la)z^n
$$
where $J_{UV,n}(\tilde\la)$ is a
rational function of $\tilde\lambda$ (respectively $q^{\tilde\lambda}$).
Also, $J_{UV,0}(\tilde\lambda)=J_{\bar U\bar V}(\la)$,
where $\la$ is the $\h^*$-part of $\tilde\la$,
where $\bar U,\bar V$ are considered as $U_q(\g)$-modules.
In particular, this shows that
$J_{UV}(\tilde\lambda)$ is invertible.

Let us also define the multicomponent fusion matrix.
Namely, let $V_1,...,V_N$ be loop representations
of $U_q(\tilde\g)$, Then define
an operator
$$
J^{1..N}(\tilde\lambda):\; V_1 \ovo...\ovo V_N \to V_1 \ovo...\ovo V_N,
$$
by
$$
J^{1..N}(\tilde\lambda):=
J^{1,2...N}_{V_1,V_2\otimes...\otimes V_N}(\tilde\la)...
J^{N-1,N}_{V_{N-1},V_N}(\tilde\la).
$$
This operator can be regarded as an
element of $\text{End}(\bar V_1\otimes...\otimes \bar V_N)
[[z_2/z_1,...,z_{N}/z_{N-1}]]$.

{\bf Remark.} It is easy to see that the matrix elements of the
operator $J^{1...N}$ are the expectation values of products of
intertwining operators (i.e. the correlation functions for the
Wess-Zumino-Witten conformal field theory, for 
$q=1$, and its q-deformation, for $q\ne 1$),
which are the main objects of discussion in \cite{EFK}. In
particular, it is known that if the representations $\bar V_i$
are irreducible then the formal series
$J^{1...N}(\tilde\la,\zeta_1,...,\zeta_{N-1})$ is convergent (for
a generic $\tilde\la$) to an analytic function of $\zeta_i$ in some
neighborhood of zero (see \cite{EFK} and references therein).
However, we will not need this fact in this paper.

\subsection{The multicomponent dynamical action}

Let $V_1,...,V_N$ be loop representations of
$U_q(\tilde\g)$.
The multicomponent shifted dynamical action $w\bullet$ of $\tilde \W^a$
on $V_1\ovo...\ovo V_N$ is defined in the same way as it was
defined in the general Kac-Moody case.

Similarly to the general Kac-Moody case,
we have the following proposition.

\begin{proposition}
The operator $J^{1...N}$ conjugates the shifted dynamical action
of $\tilde \W^a$ into its shifted multicomponent dynamical action.
That is,
$$
J^{1...N} (w\bullet)= (w\circ)J^{1...N}.
$$
\end{proposition}

\begin{proof} The proof is the same as that of
Lemma \ref{com-A} and
Lemma \ref{multi}.
\end{proof}

Define
$$
{\mathcal J}^{1...N}(\tilde\la)=
J^{1...N}(-\tilde\la-\tilde\rho+\frac{1}{2}\sum_{i=1}^N h^{(i)})
$$

\begin{corollary}
The operator ${\mathcal J}^{1...N}$ conjugates the (unshifted) dynamical action
of $\tilde \W^a$ into its multicomponent dynamical action.
That is,
$$
{\mathcal J}^{1...N} (w\diamond)= (w *){\mathcal J}^{1...N}.
$$
\end{corollary}

\begin{proof} Clear.
\end{proof}

\subsection{Trigonometric KZ equations for ${\mathcal J}^{1...N}$ ($q=1$)}

In this section we will assume that $\tilde\la=(\la,k,0)$.
Let $(,)$ be the inner product
on $\g$ that we defined before, and let
$e_\alpha,f_\alpha$ be Cartan-Weyl generators
such that $(e_\alpha,f_\alpha)=1$.

Define the Drinfeld r-matrix for $\g$
$$
r=\sum_{\alpha>0}e_\alpha\otimes f_\alpha +\frac{1}{2}
\sum_i x_i\otimes x_i
$$
(recall that $x_i$ is an orthonormal basis of $\h$).

Let $V_1,...,V_N$ be loop representations of
$U(\tilde\g)$.
Define the trigonometric KZ operators on $V_1\ovo...\ovo V_N$:
$$
\nabla_i(\tilde\la)=mkz_i\frac{\partial}{\partial z_i}
+\sum_{j\ne i}\frac{z_ir_{ji}+z_jr_{ij}}
{z_i-z_j}+\la^{(i)},
$$
where the rational functions of $z_i$
on the right hand side are expanded in a power
series with respect to $z_i/z_{i-1}$.

{\bf Remark 1.} Observe that
$mkz_i\frac{\partial}{\partial z_i}+\la^{(i)}=\tilde\la^{(i)}$
on $V_1\ovo V_2\ovo...\ovo V_N$, and that 
the operator $\frac{z_ir_{ji}+z_jr_{ij}}
{z_i-z_j}$ can be thought of as the action of 
the Drinfeld r-matrix of the {\it affine} 
algebra $\tilde\g$ in the tensor product 
of two loop representations. 

{\bf Remark 2.} We note that our trigonometric KZ operators differ
from those of \cite{TV} by a change: $\lambda\to -\lambda$,
$\kappa\to -\kappa$, $z_i\to z_i^{-1}$ (apart from the overall
factor of $m$). This is the cause of a number of sign
discrepancies between \cite{TV} and this paper.

Let $\nabla_i^0(\tilde\la)$ be the ``Cartan''
part of $\nabla_i(\tilde\la)$, i.e.
$$
\nabla_i^0(\tilde\la)=mkz_i\frac{\partial}{\partial z_i}
+\sum_{j<i}\sum_{l=1}^r x_l^{(i)}\otimes x_l^{(j)}
-\sum_{j>i}\sum_{l=1}^r x_l^{(i)}\otimes x_l^{(j)}+\la^{(i)}.
$$

\begin{thm}\label{KZ}\cite{TK,FR}
(trigonometric KZ equations)
One has
$$
\nabla_i(\tilde\la){\mathcal J}^{1...N}(\tilde\la)=
{\mathcal J}^{1...N}(\tilde\la)\nabla_i^0(\tilde\la).
$$
\end{thm}

\begin{proof} This is, after some transformations, the content of
  Theorem 3.8.1 in \cite{EFK}. This is also the
multicomponent version of the ABRR equation
for affine Lie algebras, projected to the product of
loop representations (see \cite{ES}).
\end{proof}

{\bf Remark 1.} Note that in Theorem 3.8.1 of \cite{EFK}, there
is a misprint: there should be a factor $\frac{1}{2}$ in front of
$<\mu_i,\mu_i+2\rho>$.

{\bf Remark 2.} In the KZ equation
for conformal blocks, the central charge $k$
occurs in a combination $k+h^\vee$ (see \cite{EFK}).
This shift of $k$ by $h^\vee$ is absent here
because when passing from $J$
to ${\mathcal J}$, we have performed a shift by $\tilde\rho$,
which, in particular, involves a shift of $k$ by $h^\vee$.

{\bf Remark 3.} We note that the finite-dimensional analog of Theorem \ref{KZ}
(i.e. the corresponding statement for $\g$ and not for $\tilde\g$) appears 
in \cite{TV} as formula (2) (in the case $N=2$). 
We also note that formula (2) of \cite{TV} can be generalized 
to an arbitrary Kac-Moody algebra. 

\begin{thm}\label{KZinvar} The trigonometric KZ operators
$\nabla_i(\tilde\la)$ commute with the dynamical action of $\tilde \W^a$.
\end{thm}

\begin{proof}
We have seen that the operator ${\mathcal J}^{1...N}$ conjugates
the operators $\nabla_i$ to the diagonal operators $\nabla_i^0$, and
the dynamical action of the braid group to the multicomponent
dynamical action. It is easy to see that the multicomponent
dynamical action commutes with the operators $\nabla_i^0$.
This implies the desired statement.
\end{proof}

\subsection{Trigonometric qKZ equations for ${\mathcal J}^{1...N}$ ($q\ne 1$)}

Let us now describe the generalization of the content of the
previous section to the quantum case ($q\ne 1$).

Let $V_1,...,V_N$ be loop representations of $U_q(\tilde\g)$.
Let ${\mathcal R}_{ij}(z_i/z_j)\in \text{End}(\bar V_i\otimes
\bar V_j)[[z_i/z_j]]$
be the projection of the universal R-matrix.

{\bf Remark 1.} This projection is well defined. Indeed, since $c=0$ in
$V_i$ and $V_j$
by Proposition \ref{cis0}, 
the part $q^{m(c\otimes d+d\otimes c)}$ of the universal R-matrix
disappears when it is evaluated on $V_i\otimes V_j$;
thus the R-matrix defines an element of $\text{End}(\bar V_i\otimes
\bar V_j)[[z_i/z_j]]$.

{\bf Remark 2.} We note (see \cite{FR},\cite{EFK})
that if the representations
$\bar V_1,\bar V_2$ are irreducible then the R-matrix
${\mathcal R}_{12}(z_1/z_2)\in \text{End}(\bar V_1
\otimes \bar V_2)[[z_1/z_2]]$ is not only a power series,
but is actually an analytic function (for small $z_1/z_2$), which
moreover is a product of a scalar meromorphic 
(transcendental) function on $\C$
(which is holomorphic at $0$) and an operator-valued rational function.

Let $p=q^{2mk}$.

Define the trigonometric qKZ operators
$$
\nabla_i^q(\tilde\la)=
{\mathcal R}_{i+1,i}(\frac{z_{i+1}}{z_i})...{\mathcal R}_{Ni}(\frac{z_N}{z_i})
(q^{\la})_iT_{i,p}{\mathcal R}_{i1}(\frac{z_i}{z_1})^{-1}...
{\mathcal R}_{i,i-1}(\frac{z_i}{z_{i-1}})^{-1},
$$
where $T_{i,p}z_j=z_jp^{\delta_{ij}}$.

Let $\nabla_i^{q,0}$ be the ``Cartan part'' of
$\nabla_i^q$. That is,
$$
\nabla_i^{q,0}(\tilde\la)=
q^{\sum_l x_l^{i}\otimes (\sum_{j>i}x_l^{(j)}-\sum_{j<i}x_l^{(j)})}
(q^{\la})_iT_{i,p}.
$$

\begin{thm}\label{qKZ}\cite{FR} (trigonometric qKZ equations)
One has
$$
\nabla_i^q(\tilde\la){\mathcal J}^{1...N}(\tilde\la)=
{\mathcal J}^{1...N}(\tilde\la)\nabla_i^{q,0}(\tilde\la).
$$
\end{thm}

\begin{proof} This is the main result of \cite{FR};
see also Theorem 10.3.1 in \cite{EFK}
(where the simply laced case is treated). This is also the
multicomponent version of the ABRR equation
for quantum affine algebras, projected to the product of
loop representations (see \cite{ES}).
\end{proof}

{\bf Remark.} Note that in the non-simply-laced case,
the statement of the main theorem of \cite{FR} should be corrected. Namely,
the quantity $k+h^\vee$ in the quantum KZ equations should be replaced
by $m(k+h^\vee)$, which agrees with our statements here.

\begin{thm}\label{qKZinvar} The trigonometric qKZ operators
$\nabla_i^q(\tilde\la)$ commute with the dynamical action of $\tilde \W^a$.
\end{thm}

\begin{proof} Analogous to Theorem \ref{KZinvar}.
\end{proof}

\section{The dynamical difference equations}

In this section we would like to apply the material of the previous section
to deriving the dynamical difference equations from \cite{TV}.
Before we do so, we need to establish a few auxiliary results
about evaluation representations.

\subsection{The operators $\A_i(\bold z,\la)$}

Let $V$ be a loop representation of $U_q(\tilde\g)$.  Let $D$ be
the common denominator of $(\omega_i,\omega_j^\vee)$, where
$\omega_i$ are the fundamental weights of $\g$. Let
$V^e=V\otimes_{ \C[z,z^{-1}]}\C[z^{1/D},z^{-1/D}]$. The space
$V^e$ has a natural structure of a representation of
$U_q(\tilde\g)$. We call $V^e$ the extended version of $V$.

Let
  $\pi_i\in \Pi$, $w_i\in \W^a$ be the elements
  such that $t_{\omega_i^\vee}=\pi_i w_i$.
For example, if $\omega_i^\vee=\omega^\vee(\pi)$ is the minuscule
weight corresponding to $\pi\in \Pi$, then $\pi_i=\pi,
w_i=w_{[i]}$.

Let $q=1$. For $i=1,...,r$, and any loop representation $V$ of
$\tilde\g$, consider the operators $\Pi_{i,V}$ on $V^e$ given by
the formula
$$
\Pi_{i,V}=z^{-\omega_i^\vee}(A_{w_i,V}^+)^{-1}.
$$
Also, let 
 $\hat \pi_i$ denote the automorphism of $U_q(\tilde\g)$
defined by permuting the labels of the generators according to
the permutation $\pi_i\in \text{Aut}(\Gamma^a)$. 

\begin{lemma}\label{classic} There is a unique Lie algebra automorphism
$\xi_i$ of $\tilde\g$, which satisfies the equation
$$
\xi_i(a)|_{V^e}=\Pi_{i,V}a|_{V^e}\Pi_{i,V}^{-1},\ a\in \tilde\g
$$
 in all loop
representations.  

(ii) One has
$$
\xi_i(e_j)=c_{ij}e_{\pi_i(j)},\
\xi_i(f_j)=c_{ij}^{-1}f_{\pi_i(j)}, \ \xi_i(h_j)=h_{\pi_i(j)},
\xi_i(\partial)=\partial,
$$
where $c_{ij}$ are nonzero complex numbers, and $\partial$ is a
principal grading element. 

(iii) 
There exist elements $x_i\in \tilde \h$ such
that $\xi_i=\hat\pi_i\circ \text{Ad}(e^{x_i})$.
\end{lemma}

\begin{proof} The proof is easy.
\end{proof}

So, let us define
 the Hopf algebra automorphisms $\xi_i$ of 
$U_q(\tilde\g)$ for any $q$, using the same formulas for
the action of $\xi_i$ on generators. It is easy to see that 
part (iii) of Lemma \ref{classic} is valid in the $q$-case, with the 
same elements $x_i$. 
 
{\bf Remark.} We note
that for $q\ne 1$ one no longer has
$\xi_i=\text{Ad}(z^{-\omega_i^\vee}(A_{w_i,V}^+)^{-1})$.

Now let $q\ne 1$, and let $V$ be a loop representation corresponding
to an evaluation representation $\bar V$.

\begin{proposition}
There exists a unique operator $\Pi_{i,V}^q$ on $V$, which is a
q-deformation of the classical operator $\Pi_{i,V}$, such that
$$
\xi_i(a)|_{V^e}=\Pi_{i,V}^q a|_{V^e}(\Pi_{i,V}^q)^{-1},
$$
and the determinant of $\Pi_{i,V}^q$ on $V^e/(z^{1/D}-1)$ is
independent on $q$. More specifically, the first condition defines
this operator  uniquely up to a constant, while the second
condition fixes the constant.
\end{proposition}

\begin{proof} Recall that any automorphism
of an algebra acts on the set of equivalence classes of representations 
of this algebra. All we need to show is that the representation
$V$ is stable under the automorphism $\xi_i$ for any $q$.

It follows from Drinfeld's highest weight theory of finite dimensional
representations of $U_q(\hat\g)$ that there is a unique, up to a
shift of parameter,
 finite dimensional representation
of $U_q(\hat\g)$ with the same $U_q(\g)$-character as $\bar V$
(namely, all such representations have the form $\bar V(a)$ for
some $a$). On the other hand, it is easy to check that $\xi_i$
does not change the $U_q(\g)$-character of a representation. This
implies the statement.
\end{proof}

Consider now the completed tensor product
$$
V=V_1^e\ovo
V_2^e\ovo...\ovo V_N^e:= (V_1\ovo V_2\ovo...\ovo
V_N)\otimes_{\C[z_j^{\pm 1}]} \C[z_j^{\pm 1/D}]
$$
of extended loop
representations associated to evaluation representations $\bar
V_i$. Define the operators $\Pi_{i,V}$ on $V^e$ by the formula
$$
\Pi^q_{i,V}=\Pi^q_{i,V_1}\otimes...\otimes \Pi^q_{i,V_N}.
$$

 Let $\bold z=(z_1,...,z_N)$. Consider the following operators
 on $V$:
$$
\A_i(\bold z,\tilde\la)=\Pi_{i,V}^q\A_{w_i,V}(\tilde\la).
$$

It is easy to see that the operators $\A_i(\bold z,\tilde\la)$
really depend only of the components $\la,k$ of $\tilde\la$.
Also, the parameter $k$ will be fixed in the following
discussion, so we will not write the dependence on it explicitly.
Thus we will denote $\A_i(\bold z,\tilde\la)$ simply by
$\A_i(\bold z,\la)$.

{\bf Remark.} It is easy to see that the matrix elements of the
operator $\A_i(\bold z,\la)$ are Laurent polynomials of
$z_j^{1/D}$ with coefficients in rational functions of $\la$ (or
$q^\la$).

Let $\ka=mk$. Our main result in this subsection is the following.

\begin{thm}\label{commu} The operators $\A_i(\bold z,\tilde \la)$
form a holonomic system. That is,
$$
\A_i(\bold z,\la+\ka \omega_j^\vee)\A_j(\bold z,\la)= \A_j(\bold
z,\la+\ka \omega_i^\vee)\A_i(\bold z,\la).
$$
\end{thm}

\begin{proof}
 We have (dropping $V$ from the subscripts 
and $q$ from the superscripts for brevity):
$$
\A_i(\bold z,\la+\ka \omega_j^\vee)\A_j(\bold z,\la)= \Pi_i
\A_{w_i}(t_{\omega_j^\vee}\tilde\la)\Pi_j
\A_{w_j}(\tilde\la)=\Pi_i\Pi_j \xi_j^{-1}(\A_{w_i}(
t_{\omega_j^\vee}\tilde\la)) \A_{w_j}(\tilde\la)=
$$
$$
=\Pi_i\Pi_j e^{-x_j}\hat\pi_j^{-1} (\A_{w_i}(
t_{\omega_j^\vee}\tilde\la)) e^{x_j}\A_{w_j}(\tilde\la)=
\Pi_i\Pi_j e^{-x_j+\pi_j^{-1}(w_i)(x_j)}\hat\pi_j^{-1} (\A_{w_i}(
t_{\omega_j^\vee}\la))\A_{w_j}(\tilde\la)=
$$
$$
\Pi_i\Pi_j e^{-x_j+\pi_j^{-1}(w_i)(x_j)} \A_{\pi_j^{-1}(w_i)}
(w_j\tilde\la)\A_{w_j}(\tilde\la).
$$
 Now recall that in the braid group $\tilde \W^a$ we have
$t_{\omega_i^\vee}t_{\omega_j^\vee}=t_{\omega_j^\vee}
t_{\omega_i^\vee}$, and hence $\pi_j^{-1}(w_i)w_j=
\pi_i^{-1}(w_j)w_i$ (with the length of both being
$l(w_i)+l(w_j)$ in the affine Weyl group). This implies that the
product $\A_{\pi_j^{-1}(w_i)} (w_j\tilde\la)
\A_{w_j}(\tilde\la)$ is symmetric under interchanging $i$ and $j$
(the 1-cocycle relation). Thus, the theorem is equivalent to the
statement that the expression
$G_{ij}=\Pi_i\Pi_je^{-x_j+\pi_j^{-1}(w_i)(x_j)}$ is symmetric in
$i$ and $j$. But this statement is $\tilde\lambda$-independent,
so it is sufficient to prove the theorem in the limit
$\tilde\la\to \infty$ (respectively, $q^{(\tilde\la,\alpha_l)}\to
+\infty$). We can also assume that $V$ is a single loop
representation.

Now, for $q=1$, this statement is clear since the operators $\A_i$
in the limit $\tilde\la\to \infty$ are just $z^{-\omega_i^\vee}$.
In particular, for $q=1$, conjugation by $G_{ij}$ and 
conjugation by $G_{ji}$ act in the same
way on the generators of $\tilde\g$. But since $x_i$ are
$q$-independent, this action is independent on $q$. Thus, the two
actions coincide even at $q\ne 1$. Since $V$ is an irreducible
module over $U_q(\tilde\g)$, by Schur's lemma this means that
$G_{ij}G_{ji}^{-1}=C_{ij}(q)$, where $C_{ij}(q)$ are nonzero
complex numbers, and $C_{ij}(1)=1$.

Finally, taking the determinants, we find that some power of
$C_{ij}$ is $1$, so by continuity $C_{ij}=1$ also. The theorem is
proved.
\end{proof}

\subsection{The dynamical difference equations}

Theorem \ref{commu} implies the following result, which applies
both to the classical and the quantum situation.

\begin{thm} \label{holonom} Consider the system of difference equations
$$
\varphi(\bold z,\la+\kappa \omega_i^\vee)=\A_i(\bold z,\la)
\varphi(\bold z,\la),\ i=1,...,r.
$$
with respect to a function $\varphi$ of $\lambda\in \h^*$ and
$z_1,...,z_N$ with values in $\bar V_1\otimes... \otimes \bar V_N$. Then:

(i) This system is consistent, (i.e. $\A_j(\bold z,\la+\kappa
\omega_i^\vee)\A_i(\bold z,\la)= \A_i(\bold z,\la+\kappa
\omega_j^\vee)\A_j(\bold z,\la)$ for all $i,j$).

(ii) This system commutes with the KZ (qKZ) equations:
$$
\A_i(\bold z,\la)\nabla_l^q(\la)=\nabla_l^q(\la+\kappa
\omega_i^\vee) \A_i(\bold z,\la).
$$
\end{thm}

\begin{proof} Statement (i) is exactly Theorem \ref{commu}.
Statement (ii) follows from the fact that the operator
$\Pi_{i,V}$ commutes with the KZ (qKZ) equations. 
\end{proof}

{\bf Definition.} We call this system of difference equations {\it
the dynamical difference equations} for $U_q(\tilde\g)$.

\subsection{The expression of the dynamical difference equations
via the operators $B_w^+(\la)$ in the case $q=1$}

For any quantized Kac-Moody algebra, let \linebreak
$\B^+_{w}(\la):=B^+_{w}(-\la-\rho)$. Let us write the dynamical
difference equations from the previous section in terms of the
operators $\B^+_w$ corresponding to the quantum affine algebra
$U_q(\tilde \g)$, in the case $q=1$.

\begin{proposition}\label{explform} Let $q=1$.
Then the  linear operator $\A_i(\bold z,\la)$ is defined by the
formula
$$
\A_i(\bold z,\la)=\left(\prod (z_j^{-\omega_i^\vee})^{(j)}
\right)\B^+_{w_{i}}(\tilde\la),\ i=1,...,r
$$
where $\tilde\la=(\la,k,0)$. In particular, the dynamical
difference equations for $q=1$ have the form
$$
\varphi(\bold z,\la+\kappa \omega_i^\vee)= \left(\prod
(z_j^{-\omega_i^\vee})^{(j)} \right)\B^+_{w_i}(\tilde\la)
\varphi(\bold z,\la).
$$
\end{proposition}

\begin{proof} This is immediate from the previous results.
\end{proof}

From this formula, it is seen that (as was observed already in
\cite{TV}) the difference equations corresponding to {\it
minuscule} fundamental coweights are especially simple. Namely, in
this case, $w_i\in \W\subset \W^a$, which implies that the
operator $\B^+_{w_i}(\tilde\la)$ is independent of $z_i$ and of
$\ka$.

Let us now write down an explicit formula for $\B^+_w(\tilde\la)$,
$w\in \W^a$, in a representation of the form $V_1^e\ovo
V_2^e\ovo...\ovo V_N^e$, where $V_i^e$ are extended loop
representations.

For a real root $\alpha$ of $\tilde\g$, let $\bar\alpha$ be its
$\h^*$-part, and let $m_\alpha=\alpha(d)$.

For a root $\beta$ of $\g$, let $Z_\beta$ be the operator
$z_1^{\beta^{\vee(1)}/2}...z_N^{\beta^{\vee(N)}/2}$.

\begin{proposition}\label{expl} Let $\delta: w=s_{i_1}...s_{i_l}$ be a reduced
decomposition of $w$. Let $\alpha^j$ be the corresponding roots
and $m_j:=m_{\alpha^j}$. Then on the representation $V_1^e\ovo
V_2^e\ovo...\ovo V_N^e$, one has
$$
\B_w^+(\tilde\la)=\prod_{j=1}^l\biggl(Z_{\bar\alpha^j}^{m_j}\cdot
p(-\la(h_{\bar\alpha^j})-\kappa m_j-1,
h_{\bar\alpha^j},e_{\bar\alpha^j}, e_{-\bar\alpha^j})\cdot
Z_{\bar\alpha^j}^{-m_j}\biggr).
$$
\end{proposition}

\begin{proof}
The proposition follows immediately from Proposition \ref{produ}
and the definitions.
\end{proof}

\begin{corollary} The dynamical difference equations for $q=1$ coincide
with \cite{TV}, eq. (16), after the change of variables $z_i\to
z_i^{-1}$, $\ka\to -\ka$, $\la\to -\la$.
\end{corollary}

\begin{proof} The proof is by a straightforward
comparison of the two systems.
\end{proof}

\begin{corollary}\label{tv}
The dynamical difference equations (16) in \cite{TV} are
consistent and commute with the trigonometric KZ equations (in
the form of \cite{TV}).
\end{corollary}

We note that the consistency of the dynamical difference
equations was shown in \cite{TV}, Lemma 23, using
Cherednik's theory of affinization of R-matrices. The
compatibility of the dynamical difference equations with the
trigonometric KZ equations was proved in \cite{TV} for Lie
algebras of type other than $E_8,F_4,G_2$ (Theorem 24), and was
conjectured for the remaining three types (the proof of Theorem
24 uses the existence of a minuscule fundamental coweight, and
hence fails for $E_8,F_4,G_2$). Proposition \ref{tv}
implies that this conjecture is correct.

\subsection{The case $\g={\frak sl}_n$, $q\ne 1$}

Consider the case $\g={\frak sl}_n$, $q\ne 1$. To make our
picture complete, let us calculate the operators $\Pi_{i,V}^q$ for
evaluation representations of $U_q(\hat\g)$.

First of all, it is easy to check that the operators $\Pi_{i,V}^q$ 
are consistent
with tensor product 
(i.e. they agree with morphisms mapping one evaluation 
representation into a product of others). 
On the other hand, it follows from Drinfeld's highest weight theory of 
finite dimensional representations that any evaluation representation
occurs in a tensor product of shifted vector representations
(see e.g. \cite{CP}). Thus,  
it is sufficient to compute the operators
$\Pi_{i,V}$ for the vector representation.

So let $\bar V$ be the vector representation of $U_q(\hat \g)$,
evaluated at $z=1$. Recall that this representation has a basis
$v_1,...,v_n$ in which the action of the generators is defined by
the following formulas:
$$
e_i\to E_{i,i+1}, f_i\to E_{i+1,i},
h_i\to E_{i,i}-E_{i+1,i+1},
$$
where $E_{ij}$ is the
elementary matrix, and the index $i$ is understood as
an element of $\Z/n\Z$.

The crucial properties of this formula is that these formulas
are independent on $q$, and in particular are the same as those for
$q=1$. Therefore, we get

\begin{proposition} The matrices of the operators $\Pi_{i,V}^q$ in the
basis $\lbrace v_j\rbrace$ 
are independent of $q$. They are given by the formula
$$
\Pi_{i,V}^qv_j=\gamma_{ij} v_{j+i},
$$
where 
$$
\gamma_{ij}=z^{\frac{i}{n}-1},\quad i+j>n,
$$
and 
$$
\gamma_{ij}=(-z^{1/n})^i,\quad i+j\le n.
$$
(where in the last two formulas $i,j$ are integers, not elements of $\Z/n\Z$).
\end{proposition}

For example, the matrix of the operator $\Pi_{1,V}^q$ for $n=2$ is 
$$
\Pi_{1,V}^q=\left(\matrix 0&-z^{-1/2}\\ z^{1/2}& 0\endmatrix\right).
$$


\begin{thebibliography}{WW}
\normalsize

\bibitem[AST]{AST} R.Asherova, Yu.Smirnov, V.Tolstoy,
Projection operators for simple Lie groups, Theor. Math. Phys,
v.8, issue 2, 1971.

\bibitem [ABRR] {ABRR} D.Arnaudon, E.Buffenoir, E.Ragoucy, and Ph.Roche,
{\it Universal Solutions of quantum dynamical Yang-Baxter equations},
q-alg/9712037.

\bibitem[BGG]{BGG} I.N.Bernshtein, I.M.Gelfand, S.I.Gelfand,
{\it Structure of Representations Generated by Vectors of Highest Weight},
Funct. Anal.\ Appl.\ {\bf 5} (1971), 1--8.

\bibitem[BBB]{BBB}
 Babelon, O., Bernard, D., Billey, E.,\emph{A quasi-Hopf
algebra interpretation of quantum 3-j and 6-j symbols and difference
equations}, Phys. Lett. B, \textbf{375} (1996)  89-97.


\bibitem[Ch1]{Ch1} I.Cherednik,
{\it Quantum Knizhnik-Zamolodchikov Equations and Affine Root Systems},
Commun.\ Math.\ Phys.\ {\bf 150} (1992), 109--136.

\bibitem[Ch2]{Ch2} I.Cherednik,
{\it Difference Elliptic Operators and Root Systems},
Int.\ Math.\ Res.\ Notices (1995), no.\;1, 44--59.

\bibitem[CP]{CP} Chari, V., and Pressley, A., A guide to quantum groups, 
Cambridge University Press, 1995.

\bibitem[Dr]{Dr}  Drinfeld, V.G., {\it On almost cocommutative Hopf
algebras}, Leningrad Math.J. 1(2), 1990, pp.
321--342.

\bibitem[EFK]{EFK} P.Etingof, I.Frenkel, A.Kirillov,
{\it Lectures on Representation Theory and Knizhnik-Zamolodchikov Equations},
AMS, 1998.

\bibitem[ES]{ES} P.Etingof, O.Schiffmann,
{\it Lectures on the Dynamical Young-Baxter Equations}, math.QA/9908064.

\bibitem[ES2]{ES2} P.Etingof and O.Schiffmann,
{\it Twisted traces of quantum intertwiners and quantum dynamical
  R-matrices corresponding to generalized Belavin-Drinfeld
  triples}, math.QA/0003109.

\bibitem[ESt]{ESt} P.Etingof, K.Styrkas,
{\it Algebraic integrability of Macdonald operators and
representations of quantum groups},
Comp. Math., v. 114, p.125-152, 1998.

\bibitem[EV1]{EV1} P.Etingof, A.Varchenko,
{\it Exchange Dynamical Quantum Groups},
Commun.\ Math.\ Phys.\ {\bf 205} (1999), 19--52.

\bibitem[EV2]{EV2} P.Etingof, A.Varchenko,
{\it Traces of Intertwining Operators for Quantum Groups and
Difference Equations, I}, to appear in Duke Math. J,
math.QA/9907181.

\bibitem[F]{F} G. Felder, Conformal field theory and integrable
systems associated to elliptic curves,
Proceedings of the International Congress of Mathematicians,
Z\"urich 1994, p.\ 1247--1255, Birkh\"auser, 1994;
Elliptic quantum groups, preprint hep-th/9412207,
Proceedings of the ICMP, Paris 1994.

\bibitem[FV]{FV} {\it Three formulae for eigenfunctions of integrable
  Schr\"odinger operators}, hep-th 9511120.

\bibitem[FMTV]{FMTV} G. Felder, Y. Markov, V. Tarasov, A.Varchenko,
{\it Differential Equations Compatible with KZ Equations}, math.QA/0001184.

\bibitem[FR]{FR}
Frenkel I., Reshetikhin N.,
\emph{Quantum affine algebras and holonomic
difference equations},
Commun. Math. Phys. \textbf{146} (1992), 1-60.

\bibitem[K]{K} V.Kac, {\it Infinite dimensional Lie algebras},
Cambridge University Press, 1995.

\bibitem[KoSo]{KoSo} L.Korogodsky and Y.Soibelman,
{\it Algebras of functions on quantum groups},
AMS, Providence, 1998.

\bibitem[Lu]{Lu} G.Lusztig, {\it Introduction to quantum groups},
Birkhauser, Boston, 1994.

\bibitem[RS]{RS} P.Roche, A.Szenes, Trace functionals on
  non-commutative deformations of moduli spaces of flat
  connections,
math.QA/0008149, 2000.

\bibitem[Sa]{Sa} Y.Saito, {\it PBW basis of quantized universal
enveloping algebras}, Publ. RIMS, Kyoto Univ, v.30(1994),
p.209-232.

\bibitem[TK]{TK}
A.Tsuchiya, Y. Kanie, {\it Vertex operators in conformal field
  theory on $P^1$ and monodromy representations of braid group,
Conformal field theory and solvable lattice models} (Kyoto, 1986),
 Adv. Stud. pure math, v.16, Academic press, Boston, 1988,
 pp. 297-372.

\bibitem[TV]{TV}
V.Tarasov and A.Varchenko, {\it Difference equations compatible
  with trigonometric KZ differential equations},
math.QA/0002132, 2000.

\bibitem[Zh1]{Zh1}
D.P.Zhelobenko, {\it An introduction to the theory of S-algebras
  on reductive Lie algebras,} in : Representations of infinite
  dimensional Lie groups
and Lie algebras, Gordon and Breach, 1987.

\bibitem[Zh2]{Zh2}
D.P.Zhelobenko, {\it Extremal projectors and generalized Mickelsson
algebras on reductive Lie algebras\/}, Math.\ USSR, Izv.\ {\bf 33} (1989),
85--100.

\end{thebibliography}
\end{document}